\documentclass[times]{elsarticle}
\journal{Journal of Computational Physics}
\usepackage{hyperref}
\usepackage{xcolor}
\usepackage{amsfonts,amssymb}
\usepackage{caption}
\usepackage{graphicx}%
\usepackage{physics}
\usepackage{amsthm}
\usepackage{booktabs}
\usepackage[capitalize]{cleveref}
\usepackage{enumitem}
\usepackage{algorithm}
\usepackage{algpseudocode}
\usepackage{comment}
\usepackage[format=hang,singlelinecheck=false]{subcaption}
\theoremstyle{remark}

\newcommand{\rd}{\mathrm{d}}

\newcommand{\M}{\mathcal{M}}
\newcommand{\I}{\mathcal{I}}
\newcommand{\pr}{\textbf{Pr}}

\newcommand{\true}{\mathrm{true}}
\newcommand{\R}{\mathbb{R}}

\newcommand{\vecsigma}{\pmb{\sigma}}
\DeclareMathOperator{\data}{data}

\begin{document}

\begin{frontmatter}




\title{Continuous nonlinear adaptive experimental design with gradient flow}



\author[label1]{Ruhui Jin}
\author[label1]{Qin Li}
\author[label3]{Stephen Mussmann}
\author[label4]{Stephen J. Wright}
\address[label1]{Department of Mathematics, University of Wisconsin-Madison}
\address[label3]{School of Computer Science (SCS), Georgia Institute of Technology}
\address[label4]{Department of Computer Sciences, University of Wisconsin-Madison}

\date{}
\begin{abstract}

In computational inverse problems, the optimal experimental design (OED) problem seeks the best locations in time and space at which to take measurements.
We investigate the nonlinear OED problem in the context of continuously-indexed design space for the measurements.
In contrast to traditional approaches that select experiments from a finite measurement set, a continuous design space is often a better reflection of practical experimental options, where there is considerable flexibility concerning where and when to take measurements. 
The continuously-indexed space introduces computational challenges, and we address them by employing gradient-flow and optimal transport techniques, complemented by an adaptive strategy for bi-level optimization. 
Numerical results on the Lorenz 63 system and Schr\"odinger equation demonstrate that our solver identifies good measurement times / locations and achieves improved reconstruction of unknown parameters in inverse problems.
\end{abstract}


\end{frontmatter}

\section{Introduction}

In the mathematical models of physical systems, unknown parameters in the differential equations can be inferred from carefully chosen measurements of the system. 
An optimization formulation is often used to reconstruct the unknown parameter~\cite{T05,HPUU08}. 
This formulation can be stated generically as follows:
\begin{equation}
\label{eqn: continuous loss}
\vecsigma^*:= \arg\min_{\vecsigma\in \R^d}\,  \text{Loss}[\vecsigma; \rho]\, \equiv  \int_{\Omega} \left(\M(\theta; \vecsigma) - \text{data}(\theta)
\right)^2\rho(\theta)\, \dd \theta\,,
\end{equation}
where $\vecsigma\in\mathbb{R}^d$ is the parameter to be recovered and the other elements of the formulation are  described in the rest of this paragraph. 
The notation $\theta$ indicates the design variable, selected from a space $\Omega$ of possible experimental design choices. 
Specifically, $\theta$  describes the configuration of the experiment, for instance, source and detector locations and observation times. 
Often, $\theta$ is naturally ``continuously indexed," as usually seen in regression models \cite{KW59,KH21} and medical imaging \cite{K96,LR09}. 
For instance, the time at which a measurement is taken can be {\em any} time within a certain interval, not just one of a finite number of choices.
$\rho$ denotes a probability measure over $\Omega$ that guides selection of experiments. 
Empirically, more measurements are selected for a value $\theta$ when its density $\rho(\theta)$ is larger. 
The formulation \eqref{eqn: continuous loss} incorporates a Parameter-to-Output (PtO) map $\M,$ where $\M(\theta; \vecsigma)$ gives the output at measurement index $\theta \in \Omega$ when the system parameter  is $\vecsigma$. 
Corresponding to this map is the experimental measurement $\text{data}(\theta)$ which in contrast to $\M(\theta; \vecsigma)$, may contain measurement noise.
Specifically, we assume that $\text{data}(\theta)$ is generated via ground-truth parameters $\vecsigma_{\text{true}}$ with additive random noise $\epsilon$, that is, $\text{data}(\theta) = \M(\theta;\vecsigma_{\text{true}}) +\epsilon$. 

When the design measure $\rho$ is fixed, the solution of \eqref{eqn: continuous loss} is a parameter $\vecsigma^*$ that finds the smallest weighted $L^2$ mismatch 
between the output  predicted by the model and the actual observations.
This value $\vecsigma^*$ is taken to be a good approximation to the ground-truth parameter $\vecsigma_{\text{true}}$.
Reconstruction quality depends in general on the choice of $\rho$. 
Considering that experimental data and field measurements are often expensive or time-consuming to collect, it is crucial to choose $\rho$ in the most  effective possible way.
In this paper, we develop an efficient solver to find the optimal design $\rho$ that yields high-quality reconstruction of $\vecsigma^*$ at low cost. 
This {\em optimal experimental design (OED)} problem~\cite{A04,P06,HJM24} has broad applications in clinical \cite{HPR03,BDP10}, engineering \cite{BDP10} and other natural sciences. 
A critical aspect of OED is the optimality criterion that indicates desirability of a given $\rho$, such as A-optimality and D-optimality \cite{K74}. 

A major feature of our work is that we allow the design variable $\theta$ to be continuously indexed, so that the design space $\Omega$ can represent objects such as spatial, temporal, or angular coordinates. 
Consequently, $\Omega$ is a set of uncountably many configurations. 
This formulation differs  significantly from the classical setting in which $| \Omega| <\infty$, where $\rho$ can be  represented by a finite-dimensional vector. 
Our formulation requires $\rho$ to belong to a space $\pr_2(\cdot)$, which is supported on a continuous domain $\Omega$:
\begin{equation}\label{eqn:P2}\pr_2 (\Omega) = \left\{\rho~ \Big\vert~\rho(\theta) \geq 0, \forall~\theta \in \Omega, ~\int_{\Omega} \rho(\theta) \, \dd \theta = 1, ~\int_{\Omega} | \theta|^2 \rho(\theta)\, \dd \theta < \infty \right\}\,.\end{equation}
Our optimal design problem is thus formulated as the following bilevel optimization problem, with the outer level being a minimization problem over  $\rho \in \pr_2(\Omega)$ and inner layer minimizing over $\vecsigma$:
\begin{equation}
\label{eqn: oed}
\rho^{\text{opt}} := \arg\min_{\rho\in \pr_2(\Omega)} F[\rho; \vecsigma^*[\rho]]\,,\quad\text{with}\quad \vecsigma^*[\rho]~~ \text{being~the~optimum~of~\eqref{eqn: continuous loss}}\,.
\end{equation}
Here the functional $F$ is a certain criterion that measures the quality of the experimental design.
We use such classical criteria as the A- and D-optimal designs \cite{K74}, both of which are defined in terms of the spectrum of Hessian matrix of the loss  function \eqref{eqn: continuous loss}. 
Loosely speaking, a better-conditioned Hessian leads to more stable reconstruction against measurement error, due to local strong convexity in \eqref{eqn: continuous loss}. 

The continuous OED problem~\eqref{eqn: oed} raises several important issues.
\begin{itemize}
    \item 
    The outer loop in~\eqref{eqn: oed} is an optimization posed over the probability measure space $\pr_2$ over the continuous space $\Omega$. 
    Since $\pr_2(\Omega)$ is an infinite-dimensional function space, we cannot use implementations of gradient descent over Euclidean or Sobolev spaces in a straightforward way.
    We instead propose a modified gradient-based optimization solver suited to $\pr_2$.
    \item 
    The bilevel form  of  \eqref{eqn: oed} becomes difficult to handle when the PtO map $\M$ is nonlinear. 
    In a brute-force computation, each outer-loop update requires the inner loop $\vecsigma^*$ to be obtained by solving \eqref{eqn: continuous loss}, which can be costly. 
    We develop a solver that avoids the need to perform this inner-loop optimization repeatedly.
\end{itemize}

In recent years, there has been significant progress in optimal transport \cite{FG21} and gradient flow \cite{S15}. 
Techniques for optimizing over $\pr_2$ \eqref{eqn:P2} have been developed and can be applied directly to continuous OED \eqref{eqn: oed}.
By applying the implicit function theorem, we can pass the dynamics of $\rho$ to that of $\vecsigma^*[\rho]$, eliminating the need to solve the inner-loop optimization for $\vecsigma^*$ at each iteration. 
A combination of these two strategies leads to an efficient bilevel optimization solver for the design problem \eqref{eqn: oed}.

There is a rich literature on the optimal design problem, which we sketch briefly in Section~\ref{sec:hs1}, followed by an outline of our contributions in Section~\ref{sec:hs2}.
Our method combines the classical optimal design and the novel gradient flow technique; we discuss these elements in Section~\ref{sec: prelim}. 
Section~\ref{sec: method} describes the new developments of this work.
We expand on the formulation \eqref{eqn: oed} and present first the brute-force computational method and then our approach.
Numerical tests are presented in Sections~\ref{sec: test2} and \cref{sec: test1} for the Lorenz 63 system and the Schr\"odinger model, respectively.

\subsection{Related works} \label{sec:hs1}

Optimal experimental design has a long history. 
In the last century, the area came into focus when rigorous statistical formulations for design criteria were proposed and analyzed; see the review articles and books \cite{SWMW89,P06,F13,HJM24}. 
Most settings use a finite design space $ \Omega = \{ \theta_i\}_{i=1}^m$ with corresponding discrete importances $\rho(\theta_i) $, leading to the following special case of \eqref{eqn: continuous loss}:
\begin{equation}
 \label{eqn: discrete loss}
\vecsigma^*:= \arg\min_{\vecsigma\in \R^d}\,  \sum_{i=1}^m \rho(\theta_i) \, \left(\M(\theta_i; \vecsigma) - \text{data}(\theta_i)  \right)^2\,.
\end{equation}
Theoretical foundations \cite{K74} and combinatorial design algorithms were presented in this era; see, for example, \cite{A69,W70,M00}.

Modern computing power has led to significant progress in Bayesian inverse problems and uncertainty quantification \cite{CV95,C01,CS24}.
The PtO map $\M$ is governed by ordinary or partial differential equations \cite{FH12,A21} and the Bayesian formulation provides a systematic way to quantify information. 
In the special case of {\em linear} models, the design quantity \eqref{eqn: oed} is fixed and independent of $\vecsigma$. 
For the discrete-weights version of this case, classical optimization approaches combined with sparsity control \cite{HHT08,HMLT12,APSG14} can be applied, while stochastic optimization \cite{JM21}, randomized algorithms \cite{AAS18,AS18}, and surrogate modeling can improve computational efficiency.
These techniques cannot be applied to the case of continuously indexed $\Omega$. 
(The recent paper \cite{JGLW24} described a Wasserstein gradient flow method for continuous OED.) 

Nonlinearity introduces significant new challenges in that there is no closed-form expression available to define the design criterion \cite{FL13}.
Approaches that have been proposed include bilevel constrained optimization \cite{BBKS00,HHT09}, sequential design 
\cite{HM13,HM16,FJMTR20,SA22,LBM24}, and average Bayesian formulation without prior knowledge \cite{APSG16}. 
However, few works consider the case of both nonlinear model and continuously indexed $\Omega$.  
(An exception is \cite{HM14}, which seeks the optimal choice $\theta$ from $\Omega$, rather than the probability measure $\rho \in \pr_2(\Omega)$. ) 

Active learning (AL) \cite{S09}, an important area of study in the machine learning community, is closely related to OED.
AL attempts to reduce the cost of training by choosing training data adaptively.
The approach is used in deep neural networks \cite{GIG17,SS17}
and control problems \cite{B88}. 
Experimental design provides statistical criteria that can be leveraged by AL when it chooses query data.
For instance, pool-based AL augments labeled data greedily with a new label at each iteration, leading to faster convergence to the learning goal and smaller sample complexity. 
Crucial analytical guarantees \cite{CKNS15,CN08,MRTMOG22} and empirical advances \cite{AGKK21} have been established in AL. 
Other works \cite{HHGL11,GKR10} have described Bayesian active learning and adaptive design. 

\subsection{Our contributions} \label{sec:hs2}

In this paper we develop a computational scheme for continuous nonlinear optimal design. 
To the best of our knowledge, our work is the first to address nonlinear OED by optimizing a continuous design distribution $\rho \in \pr_2(\Omega)$. 
Our bilevel approach alternates between optimization of the continuous design measure $\rho$ and optimization of the inference parameters $\vecsigma.$ 
Gradient flow is used in conjunction with the Monte-Carlo particle method, which captures the design measure $\rho$, facilitating the search for those measurements that carry maximal  information volume.  
We describe both ``brute-force" and  relaxed one-step approaches for updating $\vecsigma$, the difference being in how accurately we solve the inner problem of the bilevel formulation. 
These two strategies are specified in \cref{brute-force alg} (\cref{sec:brute-force}) and \cref{alternative alg} (\cref{sec:modified}). 
The current work can be viewed as a nonlinear extension of the paper \cite{JGLW24}, which addresses the case of continuous {\em linear} OED models.
Nonlinearity requires us to add an inner layer of optimization for $\vecsigma^*$ \eqref{eqn: continuous loss}. 

We provide computational results for 
our algorithms  on various nonlinear inverse problem models, including the Lorenz 63 dynamical system \cite{L63} (\cref{sec: test2}) and the steady-state Schr\"{o}dinger equation \cite{Z22} (\cref{sec: test1}).
Besides achieving the best design criteria, the optimal design measures reveal patterns that have  interesting correspondences to features of the physical systems. 

\section{Preliminaries} \label{sec: prelim}

We review briefly two main elements of our work: Optimal experimental design (specifically the A- and D-optimal design) and Wasserstein gradient flow.

\subsection{Optimal design for linear models}
\label{subsec: linear oed}
In the linear setting, OED is often formulated via a data matrix, with the OED goal being  to select matrix rows that are most valuable for computing the solutions.
The linear model is
\begin{equation}
\label{eqn: linear model}
\M(\theta,\vecsigma)=\mathbf{A}(\theta,:) \,  \vecsigma\,,
\end{equation}
where the ``matrix" $\bf A$ may have infinitely many (possibly uncountably many) rows, as when the index $\theta$ is continuously valued. 
The task here is to select a weighted subset of the rows in $\bf A$ that allows $\vecsigma$ to be recovered stably via linear inversion.

Consider the unknowns $\vecsigma\in\R^d$ and a design measure $\rho \in \pr_2(\Omega)$. Similar to \eqref{eqn: discrete loss}, supposing that we reweight the experiment ${\bf A}(\theta,:)$ with the density $\rho(\theta)$, the reconstruction problem~\eqref{eqn: continuous loss} becomes $\min_{\vecsigma \in \R^d} \int_{\Omega} ( {\bf A}(\theta,:)\, \vecsigma - \text{data}(\theta) )^2 \rho(\theta) \dd \theta$ has an explicit analytical solution:
\begin{equation}
\label{eqn: linear sigma}
\vecsigma^* = (\mathbf{A}^\top \mathbf{A}[\rho])^{-1} \cdot \mathbf{r}[\rho]\,,
\end{equation}
where
\begin{alignat*}{2}
\mathbf{A}^\top \mathbf{A} [\rho] & = \int_\Omega \mathbf{A}(\theta,:)^\top \mathbf{A}(\theta,:) \rho(\theta) \dd \theta && \in \mathbb{R}^{d\times d}\,,\\
\mathbf{r}[\rho] & = \int_\Omega \mathbf{A}(\theta,:)^\top\text{data}(\theta) \rho(\theta) \dd \theta&& \in\R^d.
\end{alignat*}
Since both $\mathbf{A}^\top \mathbf{A} [\rho]$ and $\mathbf{r} [\rho]$ have the form of expectations, one can employ Monte-Carlo sampling to approximate them numerically.



Uncertainties and errors in $\text{data}(\theta)$  are amplified by application of the matrix inverse $(\mathbf{A}^\top \mathbf{A}[\rho])^{-1}$. 
Recalling that we assume that $\data(\theta) = {\bf A} (\theta,:)\,\vecsigma_{\text{true}}+\epsilon$ with Gaussian noise $\epsilon \sim \mathcal{N} (0, \delta^2)$, we have that  
$\text{data}(\theta)\sim\mathcal{N}\left(\mathbf{A}(\theta,:) \,  \vecsigma_{\text{true}},\delta^2\right)$ for all $\theta$.
Moreover, the value $\vecsigma^*$ obtained from \eqref{eqn: linear sigma} is also a Gaussian random variable with a modified covariance: $\vecsigma^*\sim\mathcal{N}\left(\vecsigma_{\text{true}},\delta^2\,(\mathbf{A}^\top \mathbf{A}[\rho])^{-1}\right)$. 
To minimize the sensitivity of $\vecsigma^*$  to perturbations in data, it is crucial to choose $\rho \in \pr_2(\Omega)$ so that $({\bf A^\top A}[\rho])^{-1}$ is ``small." 
Two popular OED design criteria, known as  A-optimal and D-optimal designs, quantify ``smallness" of  $({\bf A^\top A}[\rho])^{-1}$ in different ways: They optimize the trace and determinant of the variance matrix, respectively.
We have
\begin{subequations}
\label{eqn: linear oed}
    \begin{align} \label{eqn: linear oed_a}
\text{A-optimal:} \quad \rho^A&:=\arg\min_{\rho\in\pr_2(\Omega)} F^A[\rho] \equiv \Tr\left((\mathbf{A}^\top\mathbf{A}[\rho])^{-1}\right)\,, \\
\label{eqn: linear oed_d}
  \text{D-optimal:} \quad   \rho^D &:=\arg\max_{\rho\in\pr_2(\Omega)} F^D[\rho] \equiv \log(\text{Det}(\mathbf{A}^\top\mathbf{A}[\rho]))\,.
    \end{align}
\end{subequations}
Details are discussed in \cite[Section~3.1]{JGLW24}.

In this paper, we look for the counterpart of this linear design scheme for a nonlinear map $\M$. 
When $\M$ has the linear form \eqref{eqn: linear model}, we have $\nabla_{\vecsigma} \M (\theta,\vecsigma) = {\bf A}^\top(\theta,:) \in \R^d$, and the minimization problem in \eqref{eqn: continuous loss} has Hessian $\mathbf{A}^\top \mathbf{A} [\rho]$.
For the inverse problem \eqref{eqn: continuous loss} with {\em nonlinear} $\M$, we are thus motivated to consider the Gauss-Newton Hessian approximation 
$\nabla_{\vecsigma} \M(\cdot; \vecsigma)\nabla_{\vecsigma}\, \M^\top(\cdot; \vecsigma)$, 
using this quantity in place of $\mathbf{A}^\top \mathbf{A} [\rho]$ in the OED criteria. 


\subsection{Gradient flow and particle method}
\label{sec: GF}
As seen in~\eqref{eqn: oed}, the OED problem is a constrained optimization problem over the probability measure space $\pr_2(\Omega)$. 
It calls for a numerical strategy that honors the geometry of $\pr_2(\Omega)$.

Given the optimization problem
\[
\min_{\rho\in\pr_2(\Omega)}F[\rho]\,,
\]
a gradient-based strategy demands a proper definition of the ``gradient". 
We can equip the probability measure space $\pr_2(\Omega)$ with the Wasserstein-2 metric, for which distance can be defined using the following standard definition:
\begin{equation}
\label{eqn: W2}
W_2(\mu, \nu) = \inf_{\gamma \in \Gamma (\mu, \nu)}\left(\int_{\Omega\times \Omega} \| x-y\|_2^2 \,\dd \gamma (x,y) \right)^{1/2}.
\end{equation}
The formula for  Wasserstein-2 gradient flow is thus
\begin{equation}
\label{eqn: gradient flow}
\partial_t \rho =-\nabla_{W_2} F[\rho]= \nabla_\theta \cdot \left(\rho\, \nabla_\theta \frac{\delta F[\rho]}{\delta \rho}\right)\,,
\end{equation}
where $\frac{\delta F}{\delta \rho}: \Omega \to \R$ is the Fr\'echet derivative, and $\nabla_\theta \frac{\delta F}{\delta \rho}$ is usually termed the velocity field. 
A nice property of this formula is 
that it can be represented easily by the motion of discrete particles. 
Denoting the empirical measure on selected particles $\{\theta_i\}_{i=1}^N$ by
\begin{equation}
\label{eqn:particle_presentation}
\rho=\frac{1}{N}\sum_{i=1}^N \delta_{\theta_i}\,,
\end{equation}
the McKean-Vlasov formulation states that each particle $\theta_i$ moves according to the flow: 
\begin{equation}
\label{eqn: particle v}
\dot{\theta}_i = -\nabla_\theta \frac{\delta F}{\delta \rho}(\theta_i),\, \quad  i = 1, \dots, N.
\end{equation}
The complication may arise that~\eqref{eqn:particle_presentation} is a measure-valued solution, so that \eqref{eqn: gradient flow} has to be interpreted in the weak sense; see \cite{AGS05} for a discussion of this point. 
Although the numerical computation is done via finite particle simulation,  the $W_2$ gradient flow \eqref{eqn: W2 flow} is fundamental to solving our continuous OED problem and encompasses both continuous and discrete types of probability distribution $\rho$. 

\section{Methodology for nonlinear design }
\label{sec: method}

We now describe our approaches for the nonlinear OED problem~\eqref{eqn: oed}, with the objective $F$ defined by nonlinear extensions of the familiar criteria \eqref{eqn: linear oed} from the linear OED setting.
Gradient flow (described in \cref{sec: GF}) is used for the outer-level optimization in \eqref{eqn: oed}; we provide the formulae for calculating the gradients. 
The Monte Carlo particle method is instrumental in the outer-level calculation, with particle velocities \eqref{eqn: particle v} being evaluated at the solution $\vecsigma^*[\rho]$ of the inner-level optimization of \eqref{eqn: oed}.
We describe two approaches, a ``brute-force" approach (\cref{brute-force alg}) in which the inner-level optimization problem is solved accurately, and a ``streamlined" approach (\cref{relaxed alg}) in which just one step of gradient descent is taken at the inner level, for each outer-level update. 




\subsection{Nonlinear design formulation}
\label{subsec: nonlinear oed}


The Fisher information matrix $\I[\rho; \vecsigma]$ for the mapping $\M(\theta;\vecsigma)$ is 
\begin{equation}
\label{eqn: continuous Fisher}
\begin{split}
    \I[\rho; \vecsigma] & = \int_{\Omega} \mathcal{C}(\theta; \vecsigma)\,\rho(\theta)\, \dd \theta\in\R^{d \times d} \\
    \mbox{where} \quad 
\mathcal{C}(\cdot; \vecsigma) & = \nabla_{\vecsigma} \M(\cdot; \vecsigma)\nabla_{\vecsigma}\, \M^\top(\cdot; \vecsigma) \in \R^{d \times d}\,.
\end{split}
\end{equation}
As mentioned above, the matrix $\mathcal{I}[\rho; \vecsigma]$ (sometimes also referred to in this context as the Gauss-Newton Hessian) is a natural extension of the matrix ${\bf A}^\top {\bf A} [\rho]$ that arises in the case of linear OED, which is the basis of the A-optimality and D-optimality criteria.
These criteria can be extended to the nonlinear case as follows (cf. \eqref{eqn: linear oed}):
\begin{equation}\label{eqn:FAFD}
F^A[\rho; \vecsigma] \equiv \text{Tr} \left(\I[\rho; \vecsigma] \right)^{-1}\,,\quad F^D[\rho; \vecsigma] \equiv \log(\text{Det} \left(\I[\rho; \vecsigma] \right))\,.
\end{equation}
When $\mathcal{I}[\rho;\vecsigma]$ is evaluated at the inverse problem solution $\vecsigma = \vecsigma^*[\rho]$ \eqref{eqn: continuous loss}, it is proportional to the variance of error in the solution when data is polluted by Gaussian noise. 
(The proof of this claim follows the linear case of \cref{subsec: linear oed}.)
The nonlinear experimental design problems can thus be written as follows:
\begin{subequations}
\label{eqn: A,D-optimal_star}
\begin{align}
\displaystyle \text{A-optimal:} \quad \rho^A & := \arg\min_{\rho\in \pr_2(\Omega)} F^A[\rho; \vecsigma^*[\rho]] \,,\label{eqn: A-optimal_star} \\
\displaystyle \text{D-optimal:} \quad \rho^D & := \arg\max_{\rho\in \pr_2(\Omega)} F^D[\rho; \vecsigma^*[\rho]]\,,\label{eqn: D-optimal_star}\end{align}\end{subequations}
where $\vecsigma^*[\rho]$ is obtained from the reconstruction \eqref{eqn: continuous loss}.
Because the nonlinear OED criteria~\eqref{eqn: A,D-optimal_star} depends only on the Fisher information at the global optimal point $\vecsigma^\ast[\rho]$ of the loss function, it captures the local sensitivity close to this point.

\subsection{Brute-force algorithm}
\label{sec:brute-force}

In our first solver for the nonlinear OED problem \eqref{eqn: A,D-optimal_star}, we apply the gradient flow method \eqref{eqn: gradient flow} to the outer-level problem,  yielding the formula
\begin{equation}
\label{eqn: W2 flow}
\partial_t \rho = \nabla_\theta \cdot \left( \rho\, \nabla_\theta \frac{\delta F[\rho; \vecsigma^*[\rho]]}{\delta \rho} \right)\,,
\end{equation}
where $F$ denotes $F^A$ or $F^D$ from \eqref{eqn:FAFD}, or some other criterion. 
Note that evaluation of the right-hand side requires knowledge of the  solution $\vecsigma^*[\rho]$ of the inner-level problem \eqref{eqn: continuous loss}. 
To obtain a tight approximation to this point, for each each iterate $\rho$ arising from the outer-level gradient-flow process, we run a fixed number $T'$ of gradient-descent iterations for \eqref{eqn: continuous loss}. 
We stress that \eqref{eqn: W2 flow} may be highly nonlinear, as both $\vecsigma^*[\rho]$ and the velocity field $\nabla_\theta \frac{\delta F[\rho; \vecsigma^*[\rho]]}{\delta \rho}$ depend on $\rho$ in nontrivial ways. 
For the two design objectives in \eqref{eqn:FAFD}, the velocity fields are
\begin{subequations}
\label{eqn: particle velocity_rho}
\begin{align}
\nabla_\theta\frac{\delta F^A [\rho; \vecsigma]}{\delta \rho}(\theta) & = - 2 \left(\nabla_\theta \nabla_{\vecsigma} \M^\top(\theta; \vecsigma)\right) (\mathcal{I}[\rho; \vecsigma])^{-2} \nabla_{\vecsigma} \M(\theta; \vecsigma) \in \R^{\emph{\text{dim}}(\Omega)}, \label{eqn: particle velocity_rho A}\\
\nabla_\theta\frac{\delta F^D [\rho; \vecsigma]}{\delta \rho}(\theta) & = 2\,\left(\nabla_\theta \nabla_{\vecsigma} \M^\top(\theta; \vecsigma)\right) (\mathcal{I}[\rho; \vecsigma])^{-1} \nabla_{\vecsigma} \M(\theta; \vecsigma)\in \R^{\emph{\text{dim}}(\Omega)}\,\label{eqn: particle velocity_rho D},
\end{align}
\end{subequations}
where $\mathcal{I}[\rho;\vecsigma] \in \R^{d\times d}$ is the Fisher information matrix of \eqref{eqn: continuous Fisher}. 
See~\cite[Proposition~3.1]{JGLW24} for details of this calculation in the case of linear $\M$. 

To implement gradient flow PDE~\eqref{eqn: W2 flow}, we use the particle method described in \cref{sec: GF}, restated here:
\begin{equation}\label{eqn: particle}
\rho=\frac{1}{N} \sum_{i=1}^N \delta_{\theta_i}\,,\quad\text{with}\quad\dot{\theta}_i= -\nabla_\theta \frac{\delta F[\rho; \vecsigma^*[\rho]]}{\delta \rho}(\theta_i), \;\; i=1,2,\dotsc,N.
\end{equation}
A forward-Euler time discretization for~\eqref{eqn: particle} with  step size $\Delta t$ yields 
\begin{equation}\label{eqn: particle update}
\begin{split}
\theta_i^{t+1} & = \theta_i^{t}+\Delta t\,\nabla_\theta \frac{\delta F[\rho^t; \vecsigma^*[\rho^t]]}{\delta \rho}(\theta^t_i)\,~~ \text{for}~i = 1, \dots, N, \\
\mbox{where} \;\; \rho^t & =\frac{1}{N}\sum_{i=1}^N\delta_{\theta_i^t}.
\end{split}
\end{equation}
By using the particles $\{\theta_i\}_{i=1}^N$, we can restate the velocity fields \eqref{eqn: particle velocity_rho} as
\begin{subequations}
\label{eqn: particle velocity}
\begin{align}
\nabla_\theta\frac{\delta F^A [\rho; \vecsigma]}{\delta \rho}(\theta) = - 2 \Big(\nabla_\theta \nabla_{\vecsigma} \M^\top(\theta; \vecsigma)\Big)\, \hat{\mathcal{I}}^{-2} \nabla_{\vecsigma} \M(\theta; \vecsigma)\in \R^{\emph{\text{dim}}(\Omega)}, \label{eqn: particle velocity A}\\
\nabla_\theta\frac{\delta F^D [\rho; \vecsigma]}{\delta \rho}(\theta) = 2 \Big(\nabla_\theta \nabla_{\vecsigma} \M^\top(\theta; \vecsigma) \Big)\, \hat{\mathcal{I}}^{-1} \nabla_{\vecsigma} \M(\theta; \vecsigma)\in \R^{\emph{\text{dim}}(\Omega)}\,\label{eqn: particle velocity D},
\end{align}
\end{subequations}
where $\hat{\mathcal{I}}$ is the discrete version of Fisher information \eqref{eqn: continuous Fisher}:
\[
{\hat{\mathcal{I}}}= \frac{1}{N} \sum_{i=1}^N \mathcal{C}(\theta_i; \vecsigma) \,\in \R^{d \times d}.
\]
These formulae represent a nonlinear extension of \cite[Proposition~3.4]{JGLW24}, where ${\bf A}(\theta,:)$ in \cite{JGLW24} is the gradient $\nabla_{\vecsigma} \mathcal{M} (\theta; \vecsigma)$ in the linear case. 




The complete algorithm is stated next.
The outer-level iteration appears in steps 2-3, whereas the inner-level iteration is steps 4-8.

\begin{algorithm}
\caption{Nonlinear OED Bilevel Optimization: Gradient Flow with Accurate Inner Loop}
\label{brute-force alg}
\textbf{input:} number of particles $N$; number of outer-level and inner-level iterations $T$ and $T'$, resp.; time steps for outer-level and inner-level iterations $\Delta t$ and $\Delta t'$, resp.; initial particles $\theta_1^0, \dots, \theta_N^0 \subset \Omega$; initial parameter $\vecsigma^{*,0} \in \R^d$ such that $\vecsigma^{*,0} \in \arg\min_{\vecsigma } \text{Loss}[\vecsigma; \rho^0]$.

\textbf{output:} estimates of optimal empirical measure $\rho$ (represented in the form \eqref{eqn:particle_presentation}) and corresponding inference parameter $\vecsigma^*[\rho]$. 

\begin{algorithmic}[1]
 
 \For{$t = 0:T-1$}
\vspace{1mm}

 
\State $\theta_i^{t+1} \leftarrow \theta_i^{t} + \Delta t\, \dot{\theta_i} \Big\vert_{\rho^t, ~\vecsigma^{*,t}} $, see \eqref{eqn: particle update}-\eqref{eqn: particle velocity}, with $\vecsigma^{*,t}$ replacing $\vecsigma^*[\rho^t]$, \; $i = 1, \dotsc, N$


\State $\rho^{t+1} \leftarrow \frac{1}{N} \sum_{i=1}^N \delta_{\theta_i^{t+1}}$ \Comment{(I) update measure $\rho$}

\State \textbf{initialize:} $\vecsigma^{0} \leftarrow \vecsigma^{*, t}  $\label{init line}

\For{$t' = 0:T'-1$} 
\State $\vecsigma^{t'+1} \leftarrow \vecsigma^{t'} - \Delta t' \,\nabla_{\vecsigma}\, \text{Loss} [\vecsigma^{t'}; \rho^{t+1}]$
\EndFor

\State\Return $\vecsigma^{*,t+1} \leftarrow \vecsigma^{T'}$   \Comment{(II) update parameter $\vecsigma^*$} 
\label{end line}

\EndFor
\end{algorithmic}

\Return $\rho \leftarrow  \frac{1}{N} \sum_i^N \delta_{\theta_{i=1}^T} $ and $\vecsigma^{*} \leftarrow \vecsigma^{*,T}.$
\end{algorithm}

For the inner-level iterations, the gradient of the loss function  in \eqref{eqn: continuous loss} with particle representation for $\theta$ is as follows:
\begin{equation}
\label{eqn: gradient loss}
\nabla_{\vecsigma}\, \text{Loss} [\vecsigma; \rho] \equiv \nabla_{\vecsigma}\, \text{Loss} (\vecsigma; \{ \theta_i\}) =\frac{2}{N} \sum_{i=1}^N \nabla_{\vecsigma} \mathcal{M}(\theta_i; \vecsigma) (\mathcal{M}(\theta_i; \vecsigma) - \text{data}(\theta_i)).
\end{equation}
Other optimization solvers could be applied to the inner-level problem, in place of gradient descent.


\subsection{Streamlined algorithm}
\label{sec:modified}

A  major part of the computation in \cref{brute-force alg} is the inner-loop computation of $\vecsigma^{\ast,t+1}$. 
On later iterations, when the step from $\rho^t$ to $\rho^{t+1}$ is small, we would expect that few steps of the inner loop would be required to find a sufficiently accurate approximation to $\vecsigma^{\ast,t+1}$. 
Accordingly, we describe a streamlined approach in which the inner iteration is reduced to a single step, which is an Euler-type step based on the evolution of $\vecsigma^*$, using the implicit function theorem to capture its dependence on $\rho$.


To derive the approach, we note that  $\vecsigma^\ast[\rho]$ satisfies first-order conditions for \eqref{eqn: continuous loss}, which are
\begin{equation*}
\displaystyle \nabla_{\vecsigma } \,\text{Loss} [\vecsigma^*[\rho]; \rho] = 2\int_{\Omega}\nabla_{\vecsigma}\M(\theta; \vecsigma^*)\, \left(\M(\theta; \vecsigma^*) - \text{data}(\theta))\right)\,\rho(\theta)\,\dd \theta = {\bf{0}} \in \R^d.
\end{equation*}
This equation implicitly links $\vecsigma^\ast[\rho]$ with $\rho$. 
To outline our use of the implicit function theorem \cite{KP02}, consider  $f(x,y(x))=0$, where $f$ is a smooth map.
By differentiating through with respect to $x$, we obtain $\frac{\partial f}{\partial x}  + \frac{\partial f}{\partial y} \frac{d y}{d x}=0$.
Multiplying by $\frac{dx}{dt}$ and using the chain rule, we have
$\frac{\partial f}{\partial x} \frac{\rd x}{\rd t} + \frac{\partial f}{\partial y} \frac{\rd y}{\rd t}=0$, which we rearrange to obtain $\dot{y} = -  (\frac{\partial f}{\partial y})^{-1} \frac{\partial f}{\partial x} \dot{x}$.
Setting $x = \rho$, $y = \vecsigma^\ast[\rho]$ and $f = \nabla_{\vecsigma} \text{Loss}[\vecsigma^*[\rho];\rho]$, we obtain
\begin{equation}
\label{eqn: sigma velocity_rho}
\frac{\rd \,\vecsigma^*}{\rd t} = -\left(\text{Hess}_{\vecsigma} \,\text{Loss}[\vecsigma^*; \rho]\right)^{-1}\,\int_\Omega \frac{\delta (\nabla_{\vecsigma}\, \text{Loss}[\vecsigma^*; \rho])}{\delta \rho}(\theta)\, \partial_t \rho(\theta) \dd \theta \in \R^d.
\end{equation}
Here $\text{Hess}_{\vecsigma} \,\text{Loss}[\vecsigma^*; \rho] \in \R^{d \times d}$ denotes the Hessian of the loss function \eqref{eqn: continuous loss} with input $\vecsigma = \vecsigma^*[\rho]$.
The formula \eqref{eqn: sigma velocity_rho} can be used to evolve $\vecsigma^*$ in time after each forward timestep has been taken in $\rho$. 

As in \cref{brute-force alg}, we represent the gradient flow~\eqref{eqn: W2 flow} using the particle method with \eqref{eqn: particle update}. 
The particle form of formula~\eqref{eqn: sigma velocity_rho} is 
\begin{equation}
\label{eqn: sigma velocity}
\dot{\vecsigma^*} = -\left(\text{Hess}_{\vecsigma} \,\text{Loss}(\vecsigma^*; \{ \theta_j\})\right)^{-1}\, \sum_{i=1}^N  \frac{\partial (\nabla_{\vecsigma}\, \text{Loss}(\vecsigma^*; \{\theta_j\}))}{\partial \theta_i}\, \dot{\theta_i}.
\end{equation}
Here we denote $\text{Loss}[\vecsigma^*; \rho]$ by $\text{Loss}(\vecsigma^*; \{\theta_j\})$ for the empirical measure $\rho = \frac{1}{N}\sum_j\delta_{\theta_j}$. 
The particle-based evolution equation \eqref{eqn: sigma velocity} takes the velocity of the particles $\dot{\theta}_j$ given by \eqref{eqn: particle velocity}.
The Hessian and gradient terms in \eqref{eqn: sigma velocity} can be calculated as follows:
 \begin{align}
 \label{eqn: hess loss}
 \text{Hess}_{\vecsigma} \text{Loss} (\vecsigma^*; \{\theta_j\}) & \displaystyle = \frac{2}{N} \sum_{j=1}^N \nabla_{\vecsigma} \M(\theta_j; \vecsigma^* )\nabla_{\vecsigma} \M^\top (\theta_j; \vecsigma^* )\\
 \nonumber
& \displaystyle ~~~+\frac{2}{N} \sum_{j=1}^N\text{Hess}_{\vecsigma} \M(\theta_j; \vecsigma^*) \left(\M(\theta_j; \vecsigma^*) - \text{data}(\theta_j) \right),
\end{align}
\begin{align}
\label{eqn: partial theta}
\displaystyle\frac{\partial (\nabla_{\vecsigma} \text{Loss} (\vecsigma^*;\{ \theta_j \} )}{\partial \theta_i}&  \displaystyle = \frac{2}{N} \nabla_{\vecsigma} \M(\theta_i; \vecsigma^*) (\nabla_\theta \M(\theta_i; \vecsigma^*) -\nabla_\theta \text{data}(\theta_i)))\\
\nonumber
& \displaystyle ~~~+ \frac{2}{N} \nabla_\theta \nabla_{\vecsigma} \M(\theta_i; \vecsigma^*) (\M(\theta_i; \vecsigma^*) -\text{data}(\theta_i)), \;\; i=1,2,\dotsc,N.
\end{align}

We summarize the procedures described above in~\cref{relaxed alg}.

\begin{algorithm}
\caption{Nonlinear OED Bilevel Optimization: Streamlined Gradient Flow}
\label{relaxed alg}
\textbf{input:} number of particles $N;$ number of iterations $T;$ time step $\Delta t;$

initial particles $\theta_1^0, \dots, \theta_N^0 \subset \Omega$ and parameter $\vecsigma^{*,0}$ s.t. $\vecsigma^{*,0} = \arg\min \text{Loss} [\vecsigma,\rho^0] $.

\textbf{output:} estimates of optimal empirical measure $\rho$ (represented in the form \eqref{eqn:particle_presentation}) and corresponding inference parameter $\vecsigma^*[\rho]$.

 \begin{algorithmic}[1]
 
\For{$t = 0:T-1$}
\vspace{1mm}

 \State $\theta_i^{t+1} \leftarrow \theta_i^{t} + \Delta t\, \dot{\theta_i} \Big\vert_{\rho^t, ~\vecsigma^{*,t}} $, see \eqref{eqn: particle update}-\eqref{eqn: particle velocity}, with $\vecsigma^{*,t}$ replacing $\vecsigma^*[\rho^t]$, \; $i = 1, \dotsc, N$


\State $\rho^{t+1} \leftarrow \frac{1}{N} \sum_{i=1}^N \delta_{\theta_i^{t+1}}$ \Comment{(I) update measure $\rho$}



\vspace{1mm}

\vspace{1mm}

\State $\vecsigma^{*,t+1} \leftarrow \vecsigma^{*,t} + \Delta t \, \dot{\vecsigma^*} \Big\vert_{ \rho^{t+1}, ~\vecsigma^{*,t}} $ \Comment{(II) update parameter $\vecsigma^*$} \label{line}

\EndFor

 \end{algorithmic}
 \Return{$\rho \leftarrow \frac{1}{N}\sum_{i=1}^N \delta_{\theta_i^{T}}, ~\vecsigma^* \leftarrow \vecsigma^{*,T}$}
 \label{alternative alg}
 \end{algorithm}

Replacing the inner loop of \cref{brute-force alg} by a single update step in \cref{relaxed alg}  reduces the cost per outer iteration considerably. 
The cost per outer loop in  \cref{brute-force alg} is $\mathcal{O}(dNT')$, being $T'$ evaluations of the gradient \eqref{eqn: gradient loss}. 
In \cref{alternative alg}, we account for the costs of calculating $\dot{\vecsigma^*}$ and $\dot{\theta_i}$, $i=1,2,\dotsc,N$ at each outer iteration as follows:
\begin{itemize}
\item Hessian \eqref{eqn: hess loss}: $\mathcal{O}(d^2 N)$;
\item Partial derivatives \eqref{eqn: partial theta} for all $i = 1, \dots, N$: $\mathcal{O} (dN)$;
\item Velocity field $\dot{\theta}_i$, for all $i=1,2,\dotsc,N$  \eqref{eqn: particle}, \eqref{eqn: particle velocity}: $\mathcal{O}(d^3 + d^2 N)$;
\item Computation of $\dot{\vecsigma^*}$ from \eqref{eqn: sigma velocity}: $\mathcal{O}(d^3)$. 
\end{itemize}
The overall cost is therefore $\mathcal{O}(d^3 + d^2 N)$ per outer iteration.


\section{Numerical experiment: Lorenz 63 system}
\label{sec: test2}

We demonstrate the performance of our numerical algorithms for optimal design based on the Lorenz 63 model~\cite{L63}. 
This physical system is a simplified mathematical model for atmospheric convection that is nonlinear and aperiodic. 
Though deterministic in our set up, the model solution exhibits chaotic behavior, which is why it is widely studied in many scientific domains.

Suppose we can take measurements along the trajectory of the solution to infer the parameters in the model. At which time points should we probe the model so that the collected data yields the most stable recovery? 
We apply the proposed algorithms of \cref{sec: method} to thoroughly address this experimental design question.



{\bf Model formulation and design.}
The Lorenz 63 model is three-dimensional and is defined by the following nonlinear differential equations:
\begin{equation}
\label{eqn: L63}
\left\{ 
\begin{array}{l}
\vspace{0.3em}
\displaystyle \frac{\dd x}{\dd \tau} = \alpha(y-x)\\
\vspace{0.3em}
\displaystyle \frac{\dd y}{\dd \tau} = x(\gamma-z)-y \\
\displaystyle \frac{\dd z}{\dd \tau} = xy-\beta z,
\end{array}
\right.   
\end{equation}
where $\vecsigma=(\alpha,\gamma,\beta) \in \R^3~ (d = 3)$ are the unknown parameters to be inferred, and $x,y,z$ are the state variables to be observed. 
The variables $x,y,z$  stand respectively for the convection rate and horizontal and vertical temperature variations.
We use $\tau$ to denote the time variable in the Lorenz model.
The model parameters $\alpha,\gamma,\beta$ are proportional to Prandtl number, Rayleigh number, and dimension of the physical layer, respectively. 
The ground-truth parameters are $\vecsigma_{\true}=(\alpha_\true,\gamma_\true,\beta_\true) = (10,28,8/3)$. 
The initial values of the  variables are fixed as $x(0) = 1.5$, $y(0) = -1.5$, and $z(0) = 25$. 
Numerical simulation of the model \eqref{eqn: L63} is conducted using forward Euler with step-size $\Delta \tau = 10^{-4}$. 

Next, we describe experimental design for Lorenz model. 
To reconstruct the model parameters $\vecsigma = (\alpha, \gamma, \beta)$, the PtO map $\M$ \eqref{eqn: continuous loss} is defined as follows:
\[
\M(\theta; \vecsigma)\,,\quad\text{with } \theta=(c,\tau)\in \Omega = \{x,y,z\} \times [0,3], \,\, \vecsigma = (\alpha, \gamma, \beta) \in \R^3 \, .
\]
Here the design variable $\theta = (c, \tau)$ has two layers of choices: the notion $c$ indicates which one of the three state variables $x$, $y$, or $z$ is to be measured, and $\tau$ denotes the corresponding observation time from the interval $[0,3]$. 
In other words, $\M$ captures one of the quantities $x(\tau)$ or $y(\tau)$ or $z(\tau)$  obtained from the model simulation \eqref{eqn: L63} based on the given parameters $\vecsigma = (\alpha, \gamma, \beta)$. 
The design goal therefore is to select the best time points at which to observe the state variables. Consequently, we express the design measure $\rho$ as a combination of marginal distributions $\rho_x, \rho_y, \rho_z \in \pr_2([0,3])$, all supported on the observation time interval. 
Further details about computation of Lorenz model design are given in \ref{sec: appendix}.

{\bf Algorithm setup.}
In \cref{sec: L63 D} and \cref{sec: L63 A}, respectively, we showcase the numerical results of our algorithms for D-optimal \eqref{eqn: D-optimal_star} and A-optimal \eqref{eqn: A-optimal_star} design criteria. 
Each experiment is run $20$ times, with each simulation having different particle sampling realizations from the same initial design measure $\rho^0$. 
We set the equal number of particles as observation times $\{\tau_{c,i}\}$ for each state variable $c \in \{x,y,z \} $. Thus the marginal distribution for the $x$ variable is formulated as  $\rho_x =\frac{1}{20}\sum_{i}\delta_{\tau_{x,i}}$, and similarly for $\rho_{y}$ and $\rho_z$.

\begin{figure}[!htb]
\hspace{-5em}    
\includegraphics[scale=0.42]{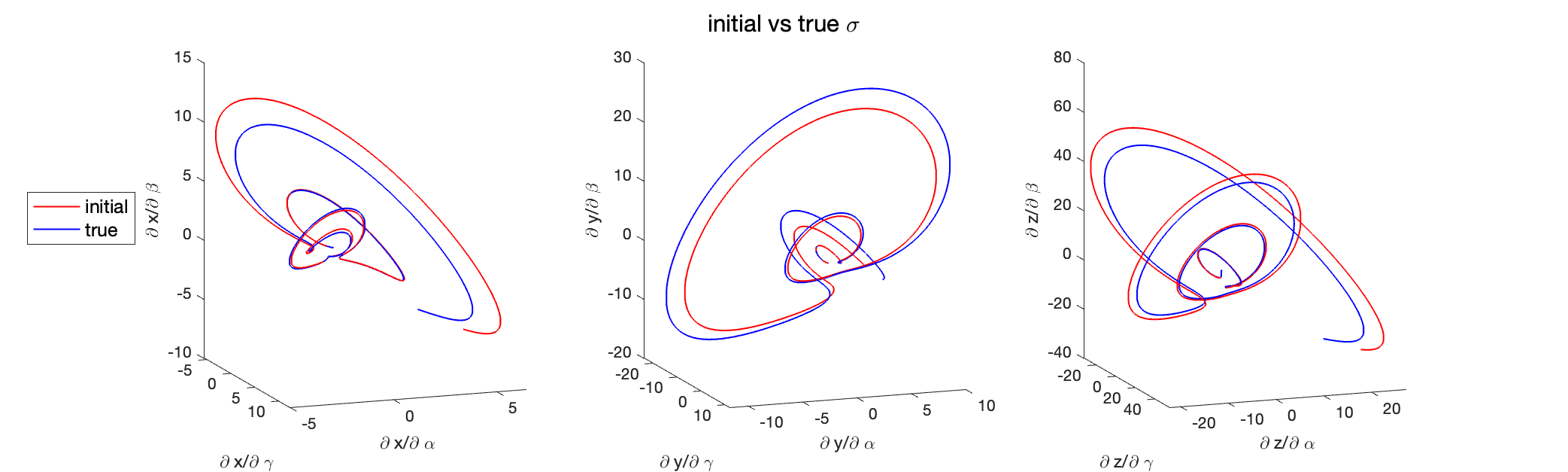}
    \caption{The comparison of $\nabla_{\vecsigma} \mathcal{M}(:, \vecsigma)$ and $\nabla_{\vecsigma} \mathcal{M}(:, \vecsigma_{\text{true}})$, where $\vecsigma \sim \vecsigma_{\true}+0.1 \, \mathcal{N}(0,{\bf I}_3)$.}
    \label{fig: initial sigma}
\end{figure}

{\bf Calculating $\vecsigma^{*,0}$.}
Given the sampled design variables $\theta = (c, \tau)$ and a starting guess of the parameters, we can realize all the measurements $\M(\theta; \vecsigma)$. 
However we must still supply the initial critical point $\vecsigma^{*,0}$ of \eqref{eqn: continuous loss}, for the given measure $\rho^0$. 
We achieve this by implementing the inner loop program in \cref{brute-force alg} for $50$ steps with step-size $10^{-5}$, with the initial guess of $\vecsigma$ drawn randomly from $\vecsigma_{\true}+0.1 \, \mathcal{N}({\bf 0},{\bf I}_3)$ for each simulation. 
Although this error gap seems to be minor, we notice an important feature of Lorenz 63 system: the model solution is highly sensitive to the parameter choice $\vecsigma = (\alpha, \gamma, \beta)$. 
As seen in~\cref{fig: initial sigma}, the perturbation following $0.1 \, \mathcal{N}({\bf 0},{\bf I}_3)$ from $\vecsigma$ to $\vecsigma_{\true}$ leads to the gradients taking on quite different values. 

{\bf Evaluating algorithm performance.}
To examine the algorithm performance, we calculate the following four quantities:
\begin{enumerate}
   \item Design objective value: $F^A$ \eqref{eqn: A-optimal_star} or $F^D$ \eqref{eqn: D-optimal_star} \label{item: F}.
    \item Parameter solution error $\|\vecsigma^* - \vecsigma_{\text{true}} \|_2$ \label{item: err}.
    \item Optimization $\text{Loss}[\vecsigma^*; \rho]$ \eqref{eqn: continuous loss}.
    \item Gradient norm $\|\nabla_{\vecsigma}\text{Loss}\|_2$.
\end{enumerate}
We will plot the evolution of these quantities, averaged over 20 simulations, in Figures~\ref{fig: L63Dmetrics}, \ref{fig: L63Dlocalmetrics}, and \ref{fig: L63Ametrics}.

\subsection{D-optimal design} 
\label{sec: L63 D}

In this section, we study the D-optimal design problem \eqref{eqn: D-optimal_star} on the Lorenz 63 model. Two sets of experiments will be conducted, each using a different initialization strategy for the design measure.

\subsubsection*{Initialization: Uniform distribution}

In the first test, we assume no prior knowledge about design preference, and so set the initialization $\rho_x^0, \rho_y^0, \rho_z^0$ to be uniformly distributed on the time window $[0,3]$. 
The total number of particles is $N = 60$, with $20$ time instances for each observable $x,y,z$. We now apply \cref{brute-force alg}, with outer-layer hyperparameters $T =50$ and $\Delta t = 10^{-5}$ and inner-layer hyperparameters $T'=20$ and $\Delta t' = 10^{-3}$.



\begin{figure}[!htb]
\hspace{-5em}    \includegraphics[scale=0.45]{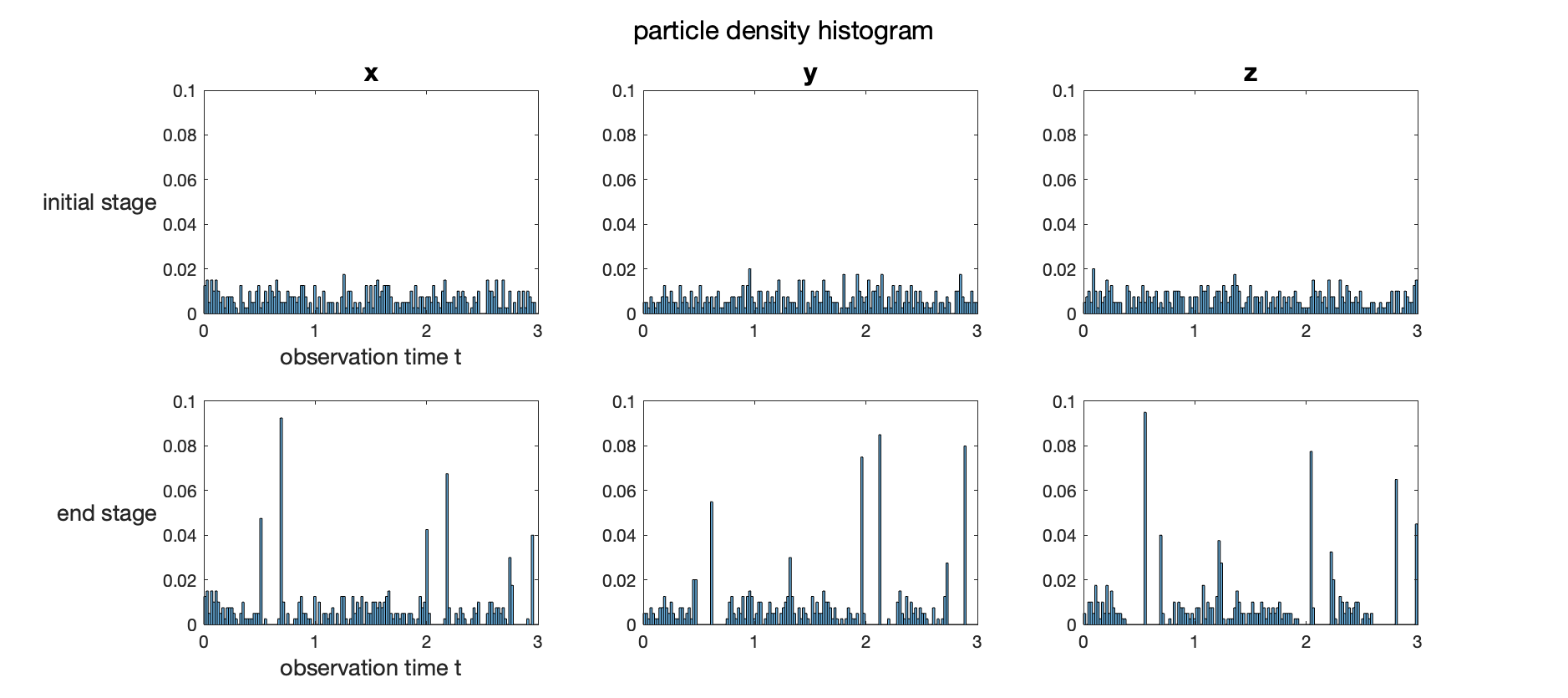}
    \caption{D-optimal design measures for $x,y,z$ observables represented as particle density histograms. The top row is the three marginal distributions of the initial phase: $\rho_x^0, \rho_y^0, \rho_z^0$, where each one is uniformly distributed over $[0,3]$. The bottom row displays the returned output of \cref{brute-force alg}.}
    \label{fig: L63Dparticle}
\end{figure}

In~\cref{fig: L63Dparticle}, we plot the marginal distributions of $\rho_x, \rho_y, \rho_z$ before and after running \cref{brute-force alg}. 
The initial design measures are uniformly distributed over $\Omega$ (top row), while in the final output of \cref{brute-force alg}, these marginal distributions eventually concentrate at specific time points (bottom row). The spikes in the density plot indicate the observation times at which measurements are crucial. 
Note that these results also suggest to observe each of the three variables at roughly the same times --- the three marginal distributions $\rho_x, \rho_y, \rho_z$  show similar spike patterns.
The histograms in \cref{fig: L63Dparticle} aggregate the particles  from all 20 simulations, meaning that we display the averaged distributions of all runs. 
Each simulation converges to a close-match design measure solution for $\rho$.

\begin{figure}[!htb]
\hspace{-5em}\includegraphics[scale=0.38]{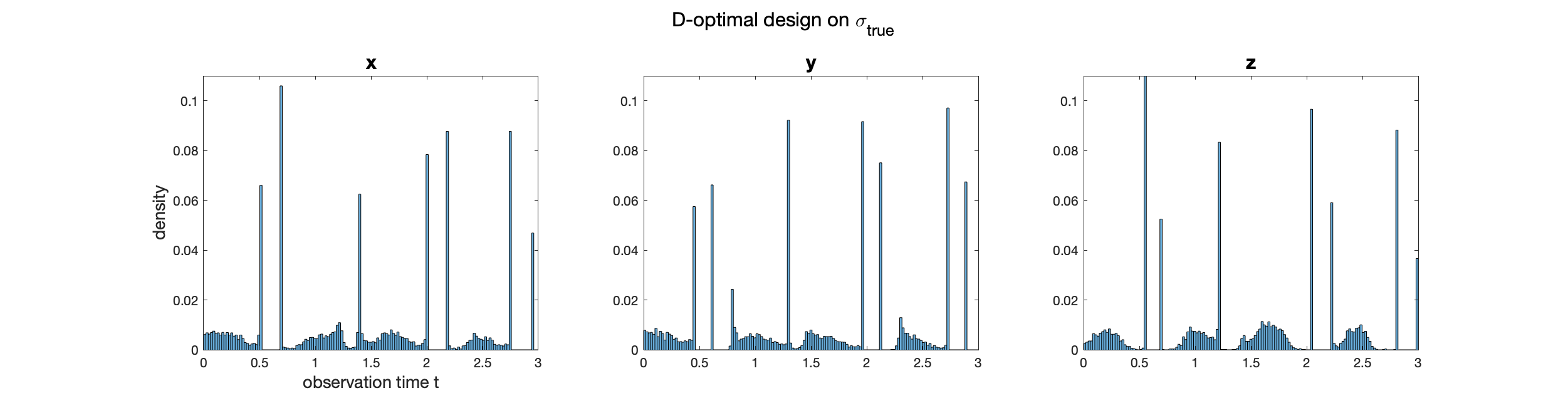}
  \caption{D-optimal benchmark design measures run on the ground-truth parameters $\vecsigma_{\text{true}}.$}
  \label{fig: L63DREF}
\end{figure}

From the known ground-truth value $\vecsigma_\true$, we can compute the reference solution for benchmarking. 
To do so, we fix $\vecsigma \equiv \vecsigma_\true$ and solve the D-optimal problem~\eqref{eqn: D-optimal_star} $\rho^{D,\text{true}}:= \arg\max_{\rho} F^D[\rho; \vecsigma_{\text{true}}]$ with large value choices of $N$ and $T$ ($N = 10000$, and $T=1000$) for the outer loop. 
(The inner loop of the two-layer optimization that updates $\vecsigma$ is not needed, since we fix $\vecsigma \equiv \vecsigma_\true$.)
In~\cref{fig: L63DREF} we showcase the density histograms of the marginal design measures $\rho_x, \rho_y, \rho_z$ returned from this one-layer design optimization process. 
The spikes exhibit the principal observation times for each state variable. 
Note the visual similarities between the benchmark result \cref{fig: L63DREF} and the adaptive design performance in the bottom row of \cref{fig: L63Dparticle}.

\begin{figure}[!htb]
\hspace{-15mm}  \subfloat[\footnotesize{D-optimal design measure highlights (circled) returned by \cref{brute-force alg} in comparison with the benchmark key observation points (all colored markers).}]{\includegraphics[scale=0.35]{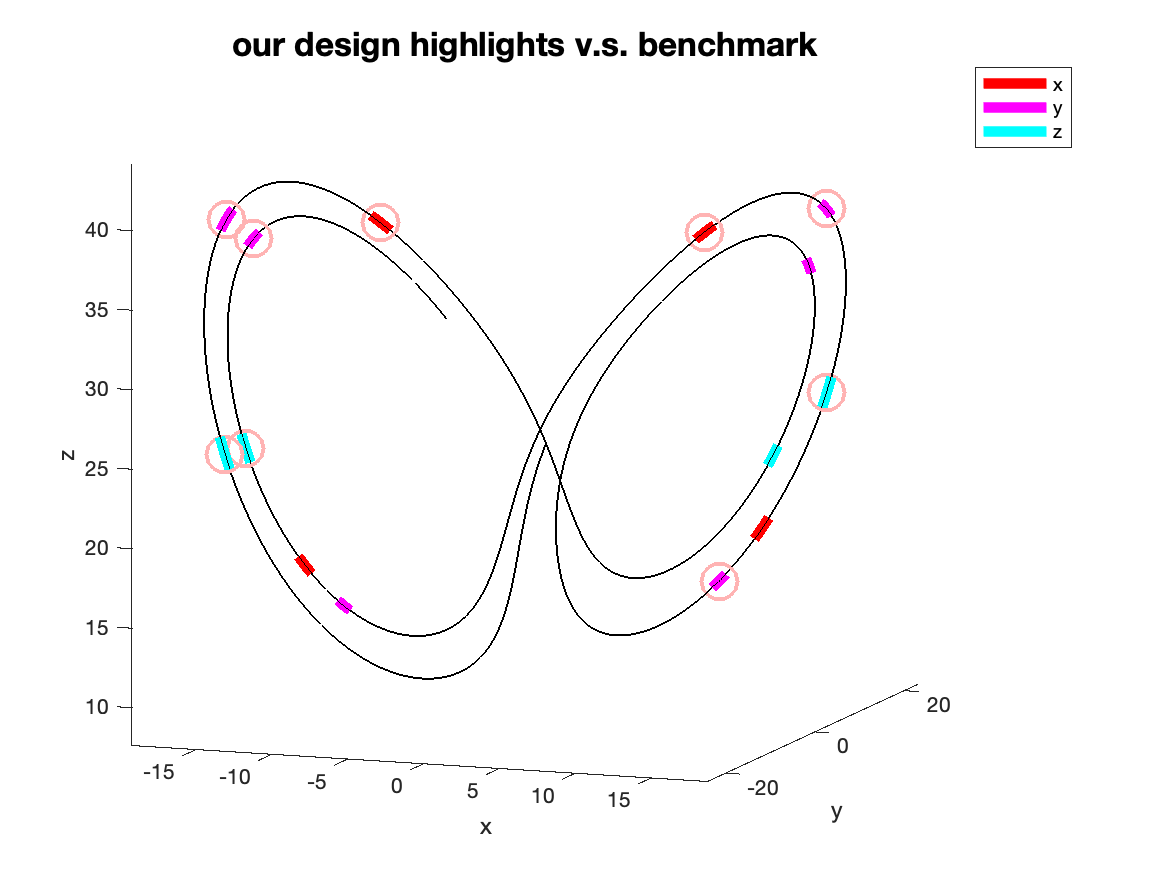}}~~
\subfloat[\footnotesize{Some importance observation times are extreme points on the model.}]{\raisebox{5ex}{\includegraphics[scale=0.35]{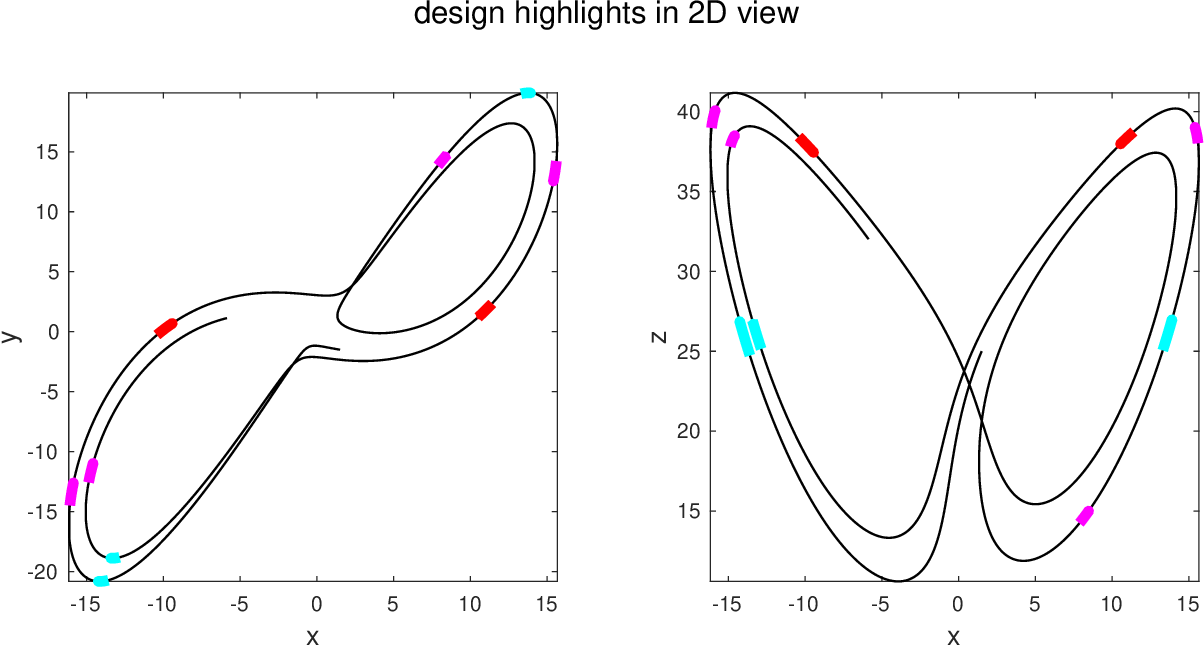}}}
    \caption{D-optimal design highlights.}
    \label{fig: L63Dmarker}
\end{figure}

For further insight, we examine the observation time points of high importance (density spikes) more closely. 
Setting the density threshold at $5\%$, we identify the time bins whose density values are above this cut-off from the histograms of \cref{fig: L63Dparticle} and \cref{fig: L63DREF}.
In~\cref{fig: L63Dmarker}, we highlight these selected time slots on the Lorenz 63 system \cref{fig: L63 model}, color-coded according to which of the three coordinates $x$, $y$, or $z$ is measured.
In part (a) of this figure, we see that the adaptive \cref{brute-force alg} result \cref{fig: L63Dparticle} recovers 9 of the 14 important observation points from benchmark computation \cref{fig: L63DREF}. 
The successful recoveries are indicated by circles.
Part (b) of the figure shows that  important measurement points tend to be located at extremal points of the trajectory of the dynamical system, in certain coordinates.

\begin{figure}[!htb]
\centering
\includegraphics[scale=0.5]{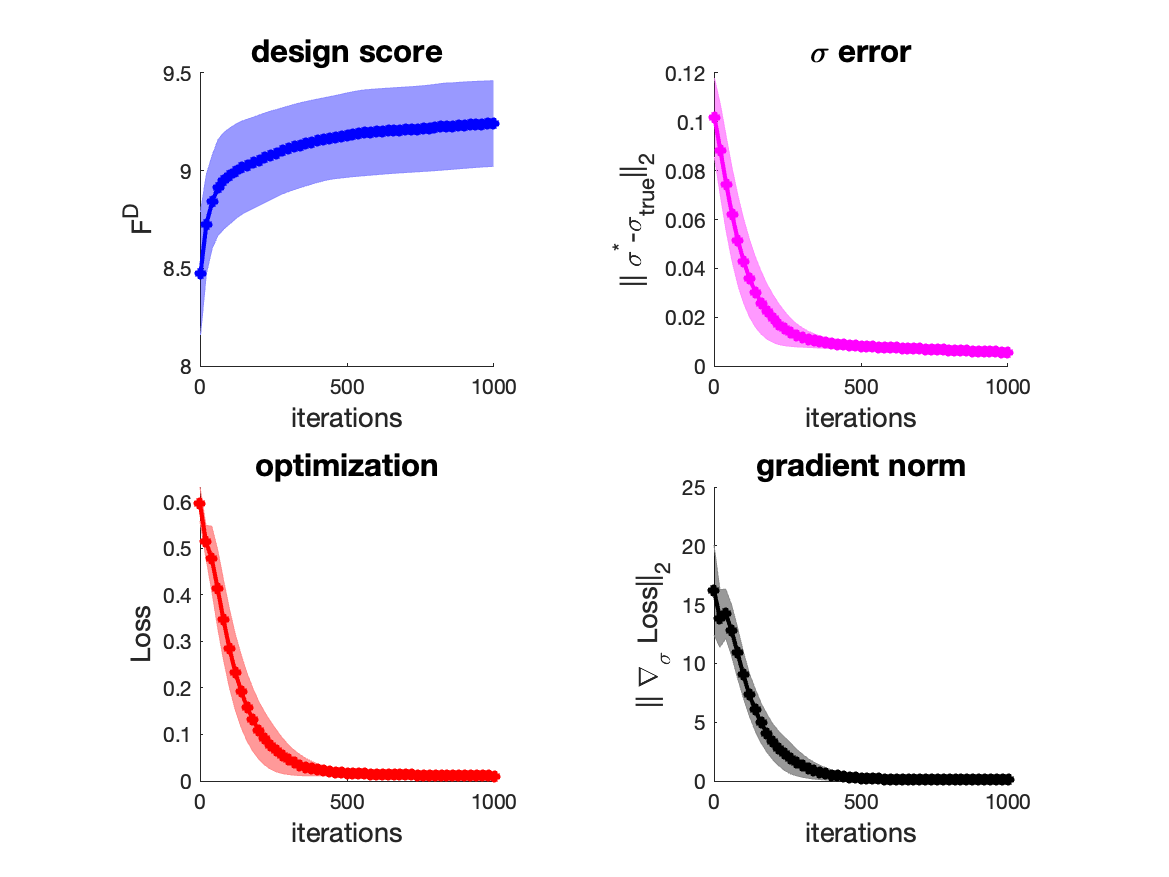}  
\vspace{-1em}\caption{Behavior of the four metrics described in the text, from uniform initial sampling $\rho^0$.}
      \label{fig: L63Dmetrics}
\end{figure}

We now quantitatively validate the performance of \cref{brute-force alg}. 
In~\cref{fig: L63Dmetrics}, we plot the evolution of four metrics representing, respectively, the D-optimal design score~\eqref{eqn: D-optimal_star}, the reconstruction error in $\vecsigma$, the loss function~\eqref{eqn: continuous loss}, and the gradient norm $\|\nabla_{\vecsigma}\text{Loss}\|_2$. 
As expected, the D-optimal design score increases to a plateau, while the  other quantities decay to values at or near zero. 
The shaded area depicts the standard deviation of $20$ independent runs. 

We now consider the quality of solution obtained by the inverse problem when observations are made at times and locations corresponding to the calculated optimal design. 
We select measurements according to two designs: our calculated optimal design (``strategic sampling") and a uniform sampling scheme. And we solve an inverse problem \eqref{eqn: continuous loss} to reconstruct the unknown parameter $\vecsigma$.
Details of these schemes are as follows:
\begin{itemize}
    \item Strategic sampling: The sampling measure is the output of~\cref{brute-force alg} shown in \cref{fig: L63Dparticle}. We only sample from the important time slots circle marked in \cref{fig: L63Dmarker}, and in each slot we sample $7$ time points randomly. Since we have identified $9$ time bins, the total sampled particles in one simulation is $N =63$.
    \item Uniform sampling: We uniformly sample $21$ points on the time window $[0,3]$ for each state variable $x,y,z$, so that the total number of sampled particles matches $N =63$.
    \end{itemize}

\begin{figure}[H]  \includegraphics[scale=0.4]{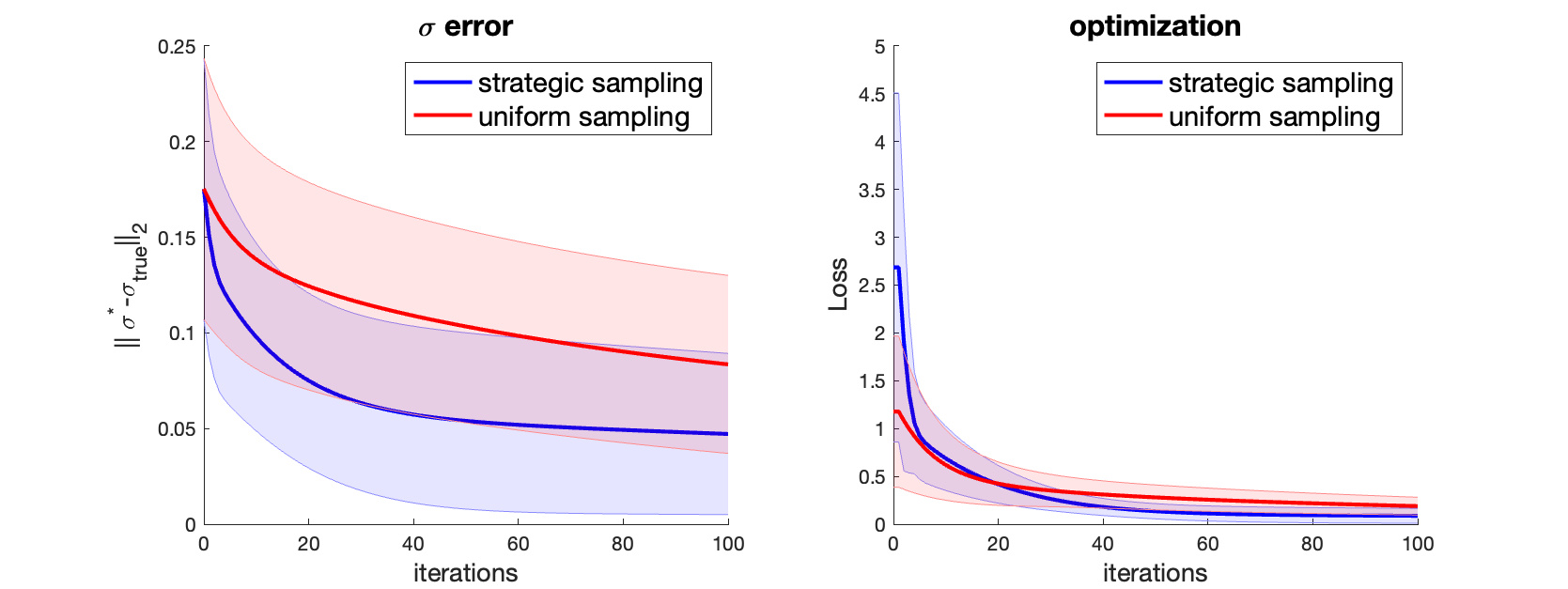}
\vspace{-1em}  \caption{Convergence comparisons between strategic sampling and a naive uniform sampling scheme.} 
     \label{fig: L63 convergence compare}
\end{figure}

The reconstruction performance comparison is shown in~\cref{fig: L63 convergence compare}. 
Strategic sampling  produces faster convergence  both for approximating the unknown model parameters in panel (a) and minimizing the reconstruction loss in (b), and also asymptotes at better final values. In both plots, the y-axis values are the averages from 20 independent simulations, where the initial parameter again follows $\vecsigma \sim \vecsigma_{\text{true}} + 0.1\, \mathcal{N} ({\bf 0},{\bf I}_3)$. The shading of each curve shows the standard deviation of all runs.

\subsubsection*{Warm-Start Initialization with Streamlined Algorithm}

Next, we study the  nonlinear design obtained from a warm-start initialization $\rho^0$. 
In this test we apply \cref{alternative alg}. 
This initial choice $\rho^0$ is inspired by the benchmark design measure shown in \cref{fig: L63DREF}. 
In particular, for each observable $x,y,z$, we select the top $6$ time bins that have the highest density values, and uniformly draw $3$ samples per bin. Thus we have $18$ samples per observable and $54$ samples in total. 

\begin{figure}[!htb]
\hspace{-6em}
\includegraphics[scale=0.45]{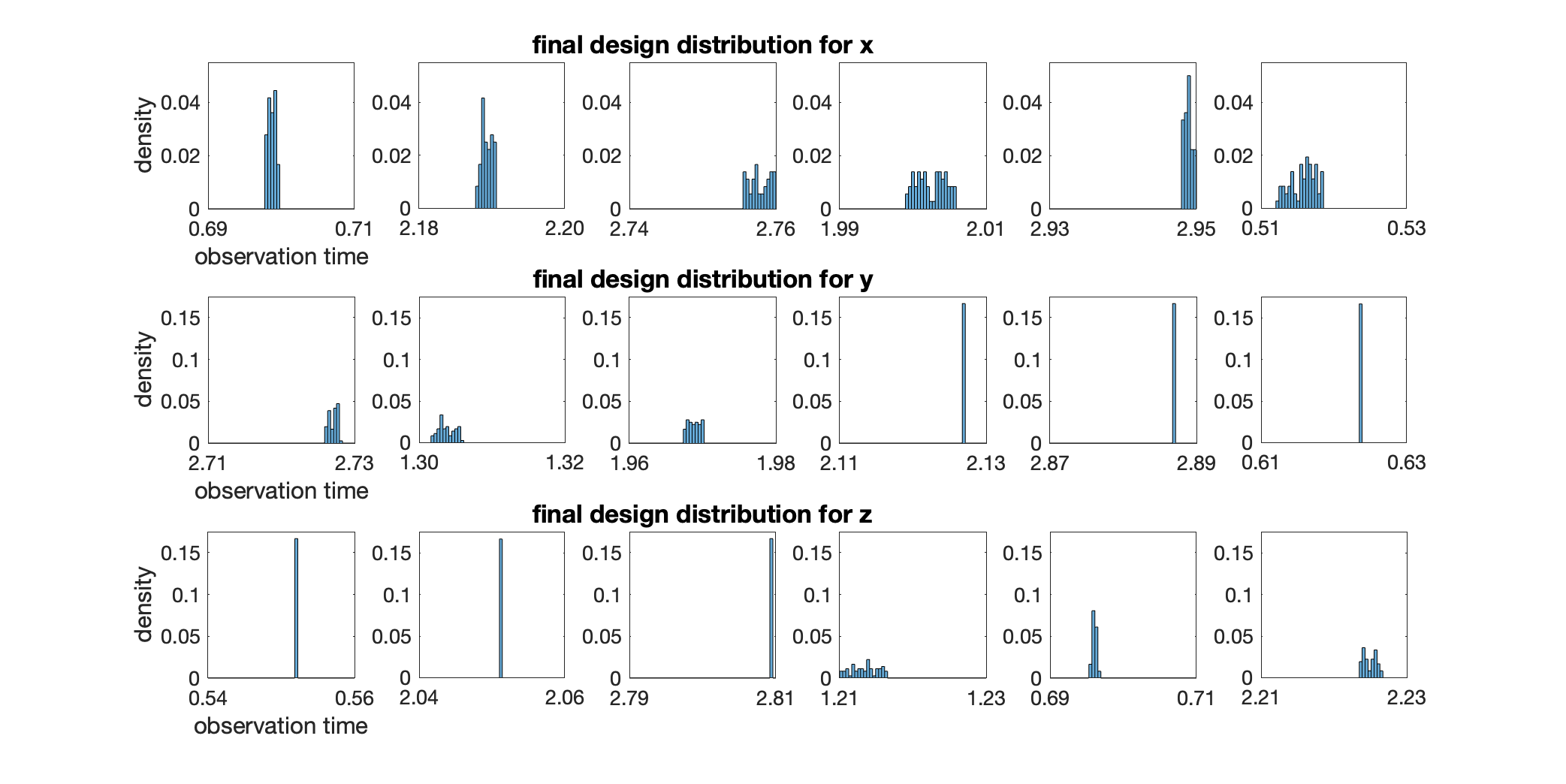}
 \vspace{-3em}   \caption{D-optimal warm-start case: Marginal distributions $\rho_x, \rho_y, \rho_z$ at the final stage of \cref{alternative alg},  where we begin by uniform sampling particles in the selected time bins for each observable $x, y, z$.}
    \label{fig: L63Dlocal}
\end{figure}

We set $T = 300$ iterations with step size $\Delta t = 0.5$ in the implementation of \cref{alternative alg}. In the initial phase, the design variables are uniformly sampled over each selected time bin (6 bins for each observable). The optimal design results are presented in~\cref{fig: L63Dlocal}: we plot the histogram for each measure $\rho_x,\rho_y, \rho_z$ in the final stage of the algorithm run restricted in these selected time bins. The design histogram is plotted counting all particle displacements from 20 independent simulations. 
The density measures to which these simulations converge are very close.

\begin{figure}[!htb]
\centering    
\includegraphics[scale=0.45]{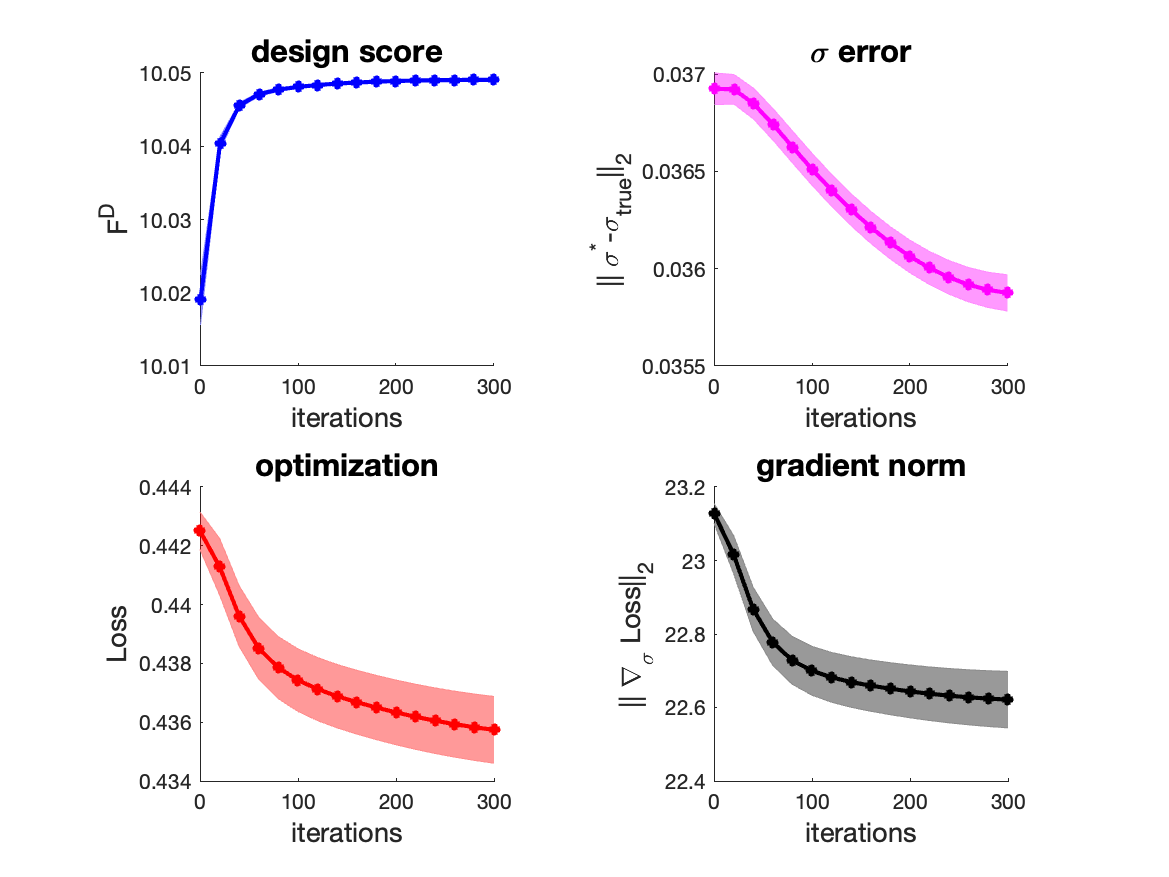}
\vspace{-1em}
\caption{Behavior of the four metrics vs iteration count when a warm start is used for $\rho^0$.}
\label{fig: L63Dlocalmetrics}
\end{figure}

We notice that the concentration continues in each time bin along the algorithm iterations. 
Although we start with samples uniformly drawn from $6$ already very narrow time windows, the output design measure has even sharper concentration, shrinking and shifting the principle observation time region further. 
The four different error metrics are plotted in~\cref{fig: L63Dlocalmetrics}, with shading showing standard deviation over 20 runs.
As expected, the D-optimal design value rises to a plateau while other error metrics decrease.



\subsection{A-optimal design}
\label{sec: L63 A}
In this section, we demonstrate the performance of proposed algorithms on A-optimal criterion \eqref{eqn: A-optimal_star}. The testing methodology follows that of D-optimal design, as described in the previous subsection.

\subsubsection*{Initialization: Uniform distribution}

\begin{figure}[!htb]
\hspace{-3em}      \includegraphics[scale=0.48]{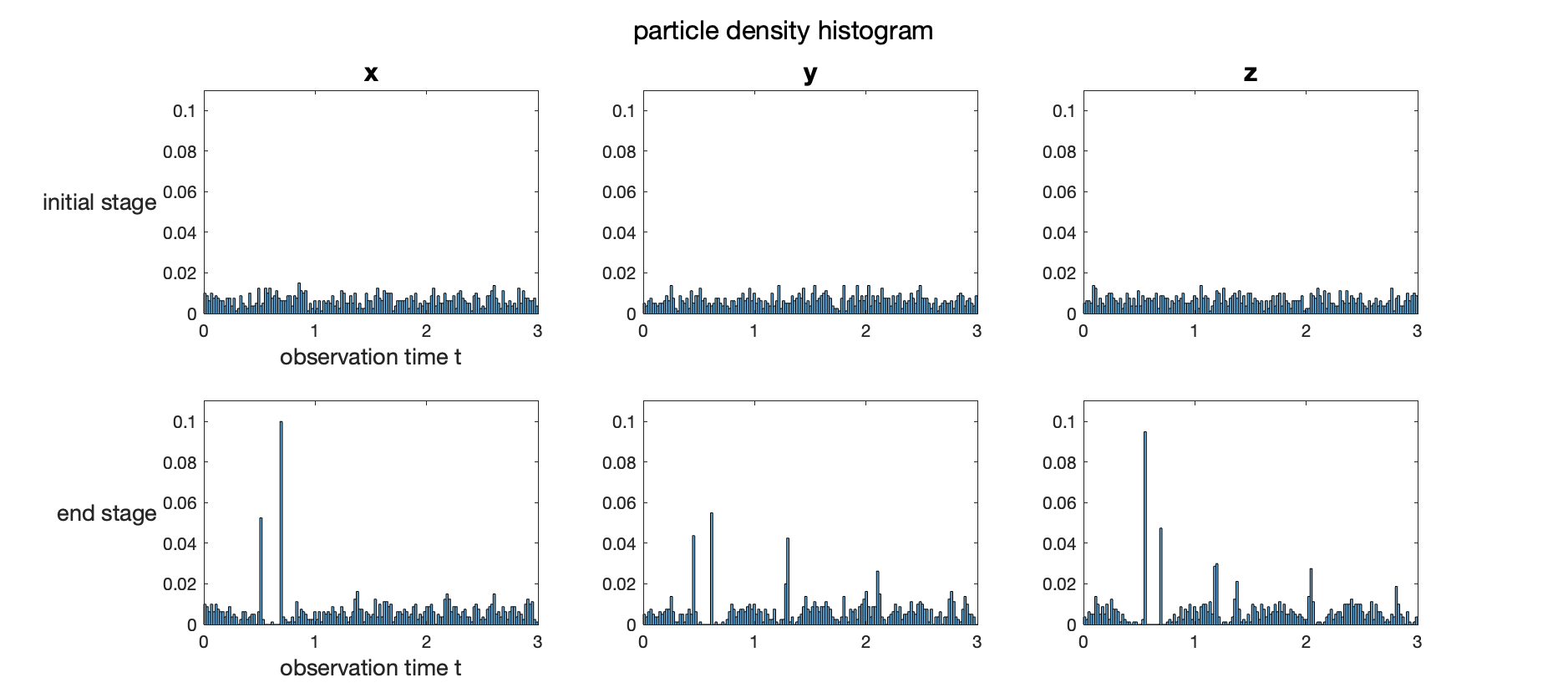}
\vspace{-2em}
    \caption{A-optimal design measures on each observable $x,y,z.$ The top row is the initial uniform distribution. The bottom row is the \cref{brute-force alg}-optimized distribution.}
    \label{fig: L63Aparticle}
\end{figure}

As the initial design measure $\rho^0$ is in the uniform pattern, we begin by sampling $20$ time points uniformly on the interval $[0,3]$ for each state variable $x,y,z.$ We then run~\cref{brute-force alg}, setting $T=50$, $\Delta t = 10^{-1}$ for the outer loop and $T' = 20$, $\Delta t' = 10^{-5}$ for the inner loop. The A-optimal design outcome is featured in \cref{fig: L63Aparticle} as the averaged density result from 20 simulations, where each individual run results in almost identical measure.
We find that the observation times concentrate in a few time windows, as in  the D-optimal result shows in \cref{fig: L63Dparticle}. 
However, there are marked differences in the time locations of the spikes between  the A-optimal and D-optimal designs (\cref{fig: L63Aparticle} and \cref{fig: L63Dparticle}, respectively).

\begin{figure}
\hspace{-2em}
\includegraphics[scale=0.35]{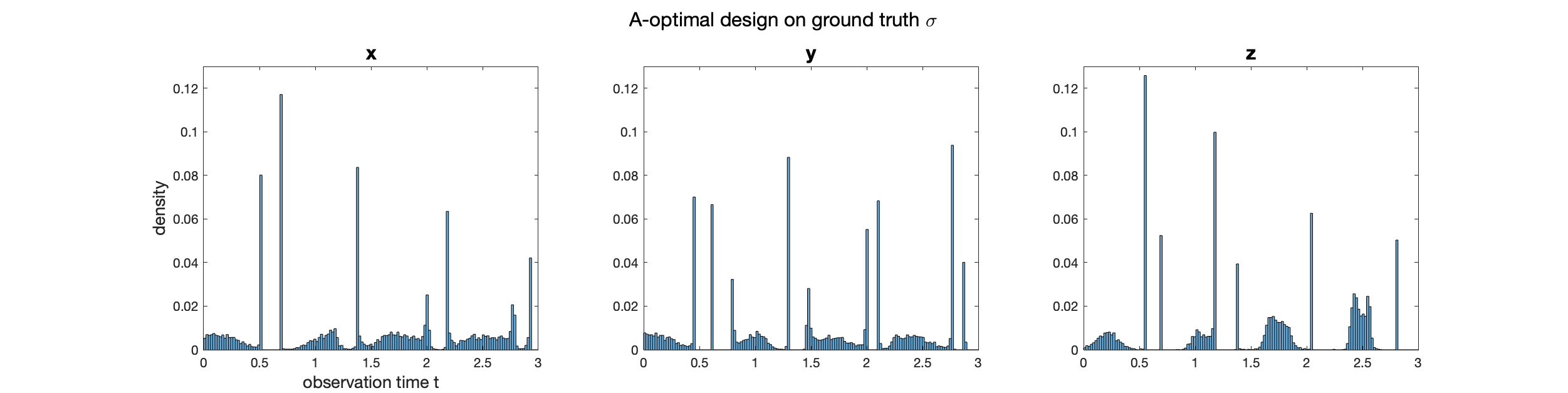}
\vspace{-1.5em}
\caption{A-optimal benchmark design measures run on the ground-truth $\vecsigma_\true$.}
\label{fig: L63Abenchmark}
\end{figure}

In~\cref{fig: L63Abenchmark} we show the benchmark A-optimal design solution by setting $\vecsigma \equiv \vecsigma_\true$ and solving   $\rho^{A, \text{true}}:= \arg\min_\rho F^A[\rho; \vecsigma_{\text{true}}]$. 
To do so, we implement only the outer loop update on $\rho$ with iterations $T = 500$ and $N = 10^4$ samples.  
Note that the benchmark density spikes of A-optimal (\cref{fig: L63Abenchmark}) are located in similarly to the D-optimal benchmark result in \cref{fig: L63DREF} (a), indicating that A-optimal and D-optimal designs are similar on the Lorenz model.

We highlight the concentrations from \cref{fig: L63Aparticle} and \cref{fig: L63Abenchmark}, and mark these temporal slots correspondingly on the Lorenz trajectory. Once again we identify the important observation times by using $5\%$ as the cut-off threshold. The benchmark design solution recognizes $10$ time windows the trajectory of Lorenz 63, and our \cref{brute-force alg} recovers $7$ of them (circled).

\begin{figure}[!htb]
\centering
\includegraphics[scale=0.35]{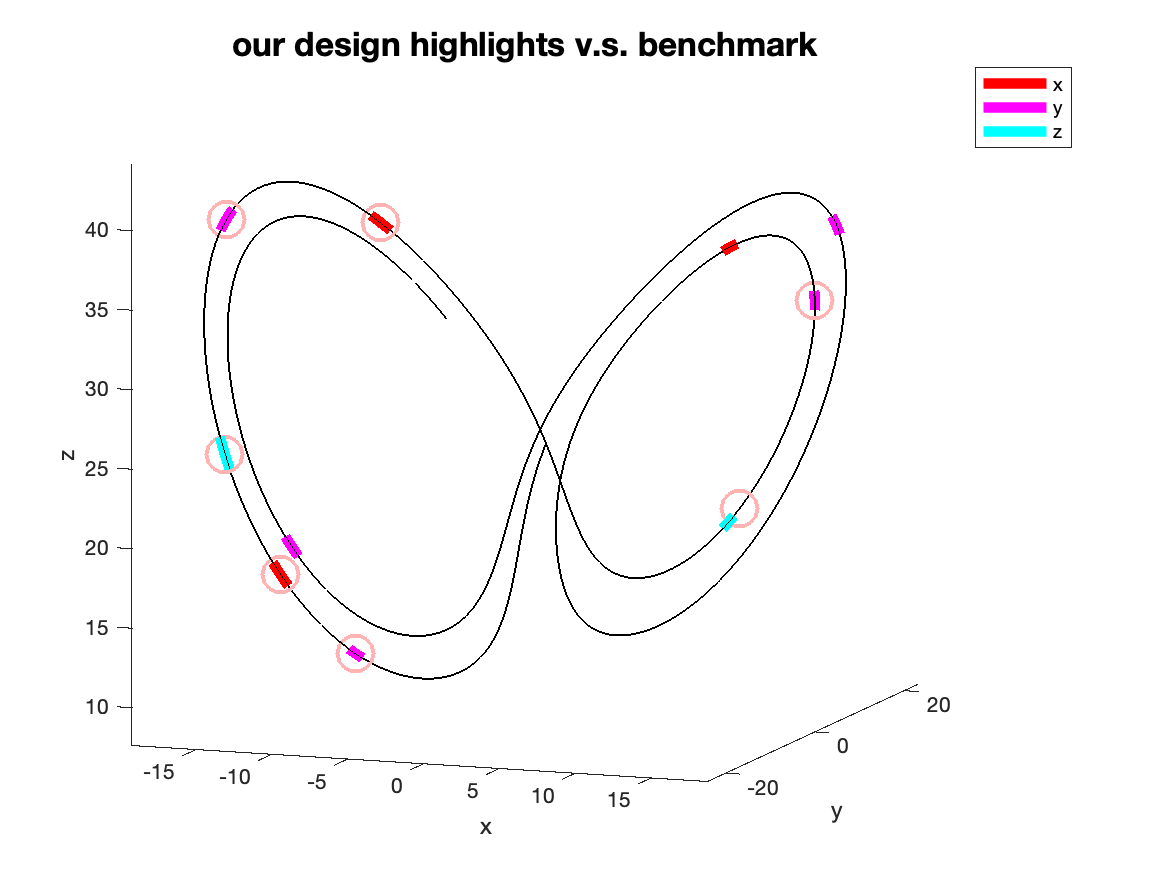}
\vspace{-1em}
    \caption{A-optimal design highlights by \cref{brute-force alg} (in red circles) recover 7 out of the total 10 benchmark important spots (all color marked slots). The three marking colors distinguish the results for $x,y,z$ observables respectively. } 
    \label{fig: L63Amarkers}
\end{figure}

\begin{figure}[!htb]
\centering
\includegraphics[scale=0.45]{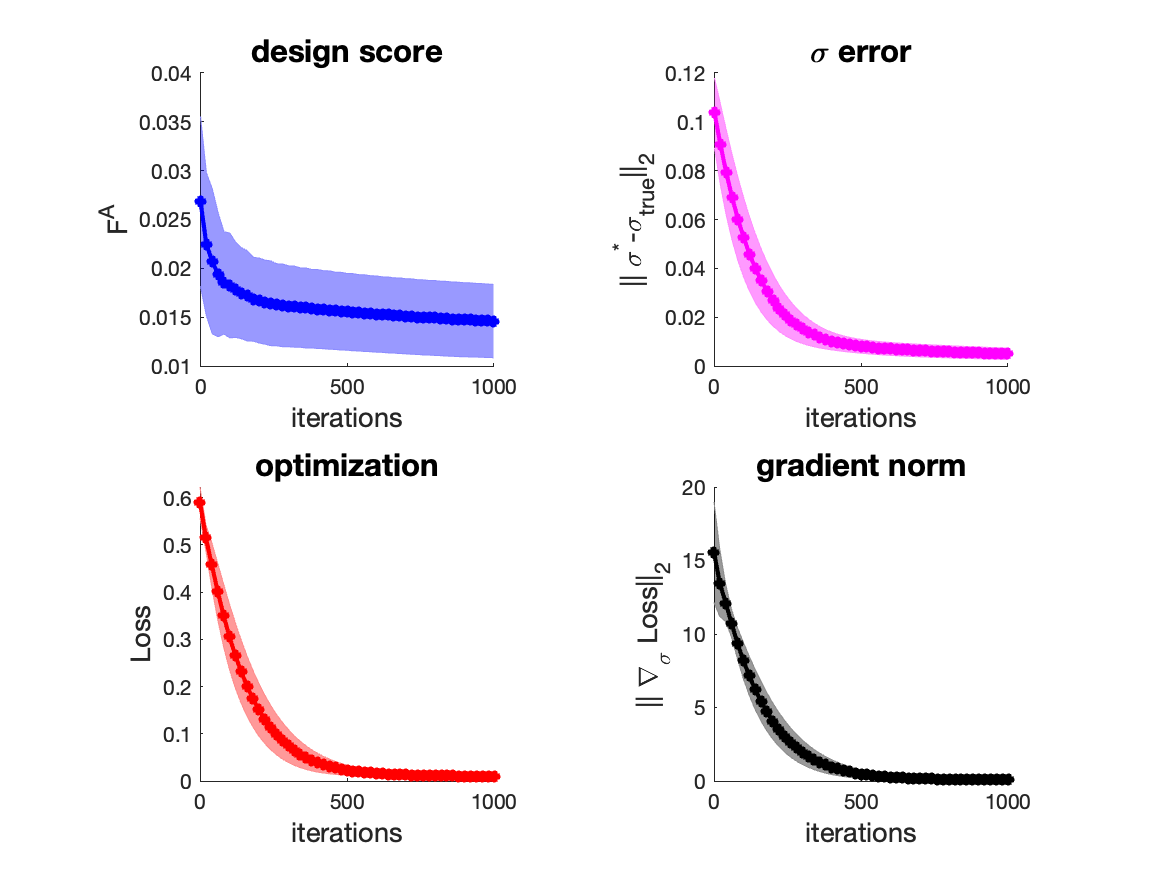}
\vspace{-1em}
\caption{The four testing metrics when the initial design measure $\rho^0$ is uniform distribution. }
\label{fig: L63Ametrics}
\end{figure}

\clearpage
The evaluation metrics for \cref{brute-force alg} are plotted in \cref{fig: L63Ametrics}. 
The A-optimal design score \eqref{eqn: A-optimal_star} decreases to a plateau, while the error in $\vecsigma$, the loss function~\eqref{eqn: continuous loss}, and gradient norm $\|\nabla_{\vecsigma}\text{Loss}\|_2$ all decay to values at or near zero.

\subsubsection*{Warm-start Initialization with Streamlined Algorithm}

Finally, we  consider a warm start for the A-optimal design measure and apply \cref{alternative alg}. 
To initialize  the algorithm, particles are sampled from only important time windows. 
In particular, we choose the top 6 time bins for each observable $x,y,z$ according to the benchmark result showed in \cref{fig: L63Abenchmark}. We sample 3 particles for each time slot, giving a total of $N =18$ samples.

\begin{figure}[!htb]
\centering
\includegraphics[scale=0.35]{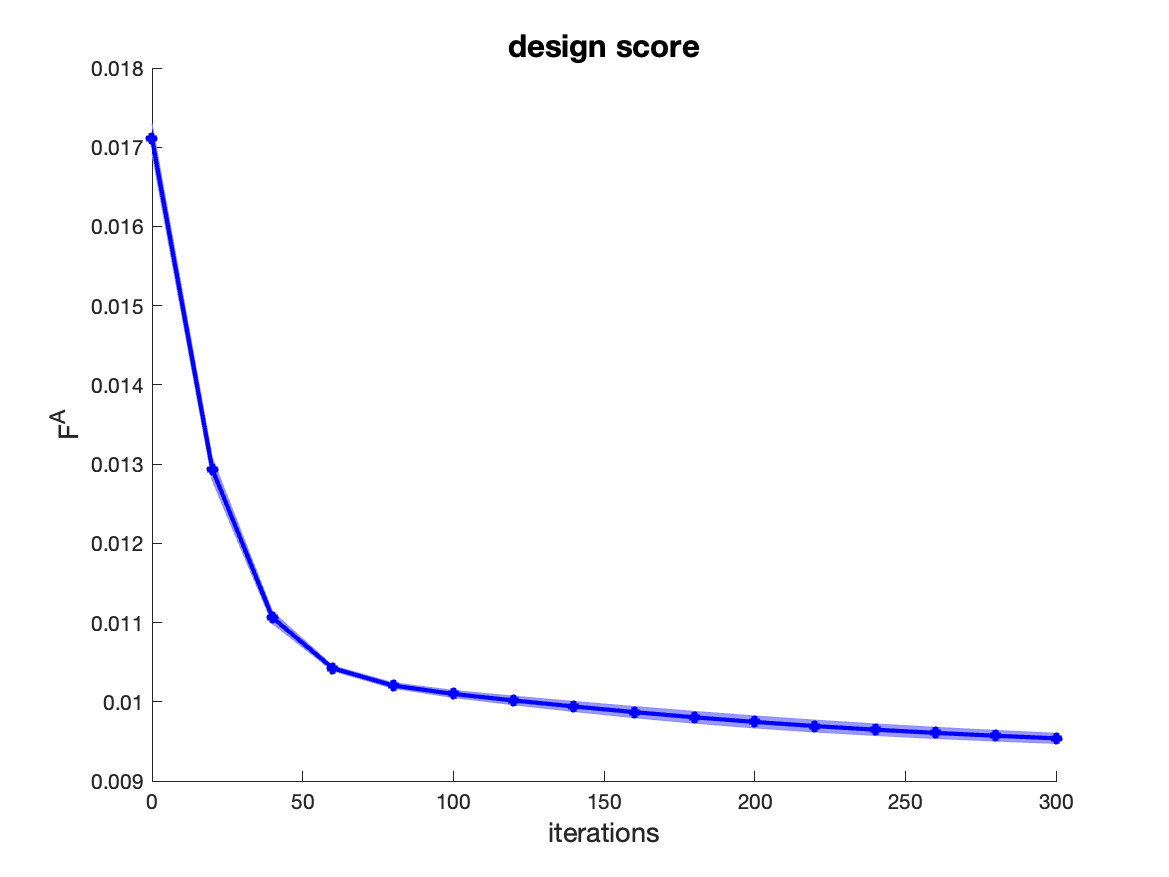}\caption{A-optimal design objective \eqref{eqn: A-optimal_star} decreases along \cref{alternative alg} iterations when the initialization $\rho^0$ is set as a warm start.}
\label{fig: L63Alocalmetrics}
\end{figure}

\begin{figure}[H]
\hspace{-2em}
\subfloat{\includegraphics[scale=0.48]{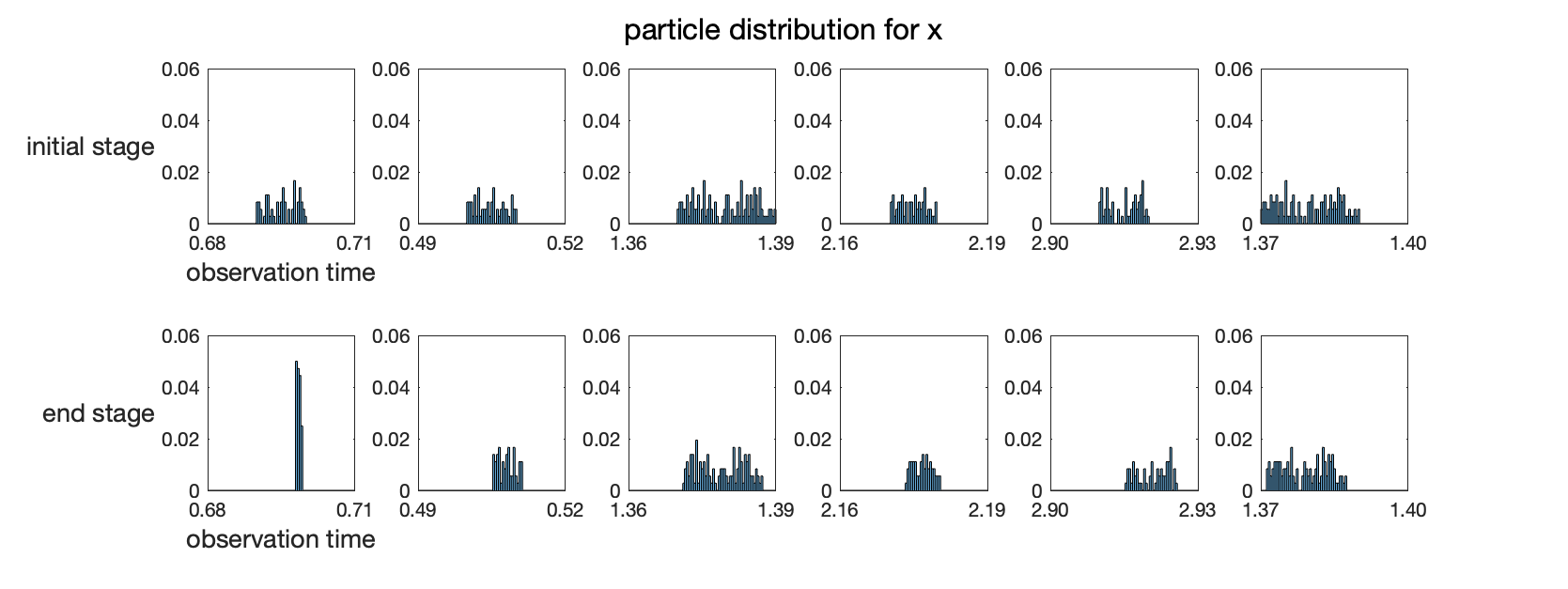}}\\
\hspace*{-2em} \subfloat{\includegraphics[scale=0.48]{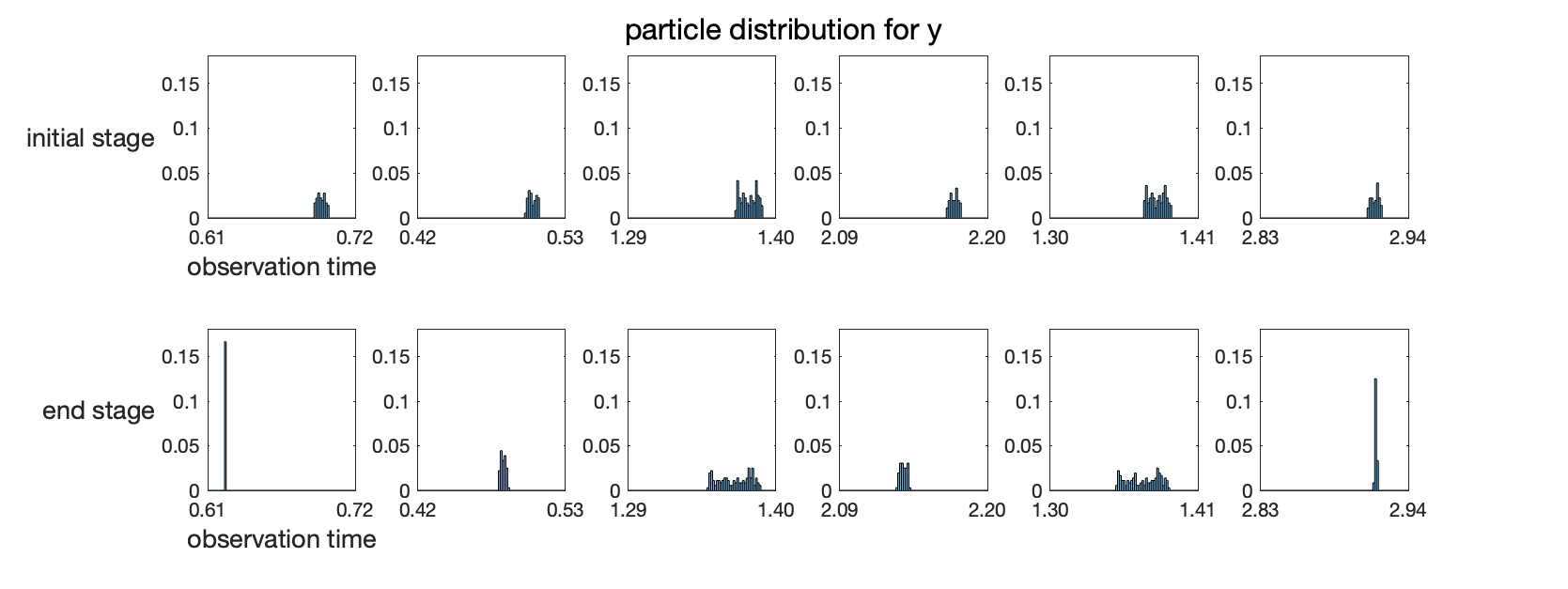}}\\
\hspace*{-2em}\subfloat{\includegraphics[scale=0.48]{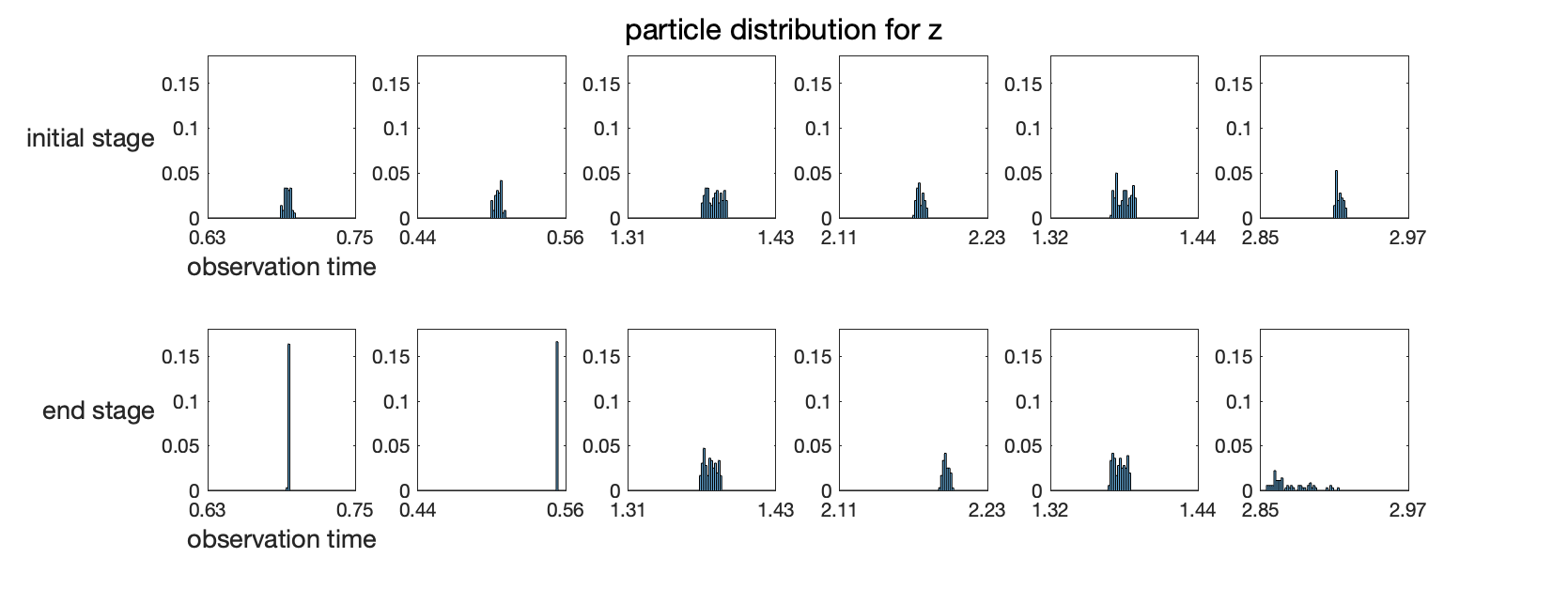}}
\vspace{-1em}
    \caption{A-optimal design warm-start case: we present the marginal distributions $\rho_x, \rho_y, \rho_z$ at the initial and final stages of \cref{alternative alg}, zooming in 6 important time bins for each observable $x,y,z$.}
    \label{fig: L63Alocal}
\end{figure}

We proceed with the iterations of \cref{alternative alg} by setting $T = 300$, $\Delta t = 50$ \footnote{The particle velocities are already very weak, so we can choose a large  step size.}. 
The marginal optimal design measures $\rho_x, \rho_y, \rho_z$ returned by the algorithm are shown in~\cref{fig: L63Alocal} as the particle density from all 20 algorithm runs. 
By comparing the initial and the final stages of the design distributions, we see that most particles remain roughly within their initial windows, and in some other bins, sometimes concentrating more sharply. 
Note too that the selected time bins for each observable $x,y,z$ are located in similar time intervals, which suggests co-dependence between the state variables in the Lorenz 63 model.
The A-optimal objective score \eqref{eqn: A-optimal_star} is shown in~\cref{fig: L63Alocalmetrics}. The value for $F^A$ drops only slightly, implying that the initial $\rho^0$ is already close to the optimum $\rho^A$, making the effect of OED update less significant. 

\vspace{-1em}

\clearpage 
\section{Numerical experiment: Schr\"{o}dinger equation}
\label{sec: test1}

We turn next to a simplified Schr\"{o}dinger equation model~\cite{G18}, for which we use~\cref{alternative alg} to identify important measurements. 
Unlike the Lorenz system, which evolves in time, the simplified Schr\"{o}dinger equation describes a steady-state solution, with measurements taken {\em in space}. 
In experimental design, we are tasked to find the measurement locations that contain the best information to recover the unknown potential.

{\bf Model formulation and design.}
The steady state Schr\"{o}dinger model with zero boundary condition is
\begin{equation}
\label{eqn: schrodinger}
\left\{ 
\begin{array}{l}
u'' + \sigma\, u = S(x)\,, \quad x \in [0,1]  \\
u \Big\vert_{x = 0, 1} = 0\,,
\end{array}
\right.
\end{equation}
where we set the domain $x\in[0,1]$. 
Different sources $S(x)$ are injected into the model and the corresponding model solutions $u(x)$ are taken to infer the unknown potential function $\sigma \in L^2([0,1])$. 
In our numerical experiment, we set the ground-truth potential $\sigma_\true$ to be a Gaussian mixture; details are given below.

For the purpose of experimental design, we can exploit the form of the source $S(x)$ and the detector location. Specifically we can inject source at a location $\theta^1 \in [0,1]$ by setting $S(x) = \delta_{\theta^1}$, and measure the solution of the resulting model \eqref{eqn: schrodinger}, denoted by $u_{\theta^1}$, at location $\theta^2 \in [0,1]$. 
The design variable is therefore $\theta = (\theta^1, \theta^2) \in\Omega = [0,1]^2$ and the simulated measurement is in the form:
\[
\M(\theta;\sigma ) = u_{\theta^1}(\theta^2) \in \R\,.
\]
Using a classical calculus-of-variable argument, one can compute the Fr\'echet derivative of $\M$ with respect to $\sigma$ to be
\[
\frac{\delta \M(\theta;\sigma )}{\delta\sigma } = u_{\theta^1} \cdot v_{\theta^2}\,,
\]
where $u_{\theta^1}$ and $v_{\theta^2}$ solve the forward and adjoint models, which are obtained from the following equations:
\[
\text{forward}: ~~\left\{ 
\begin{array}{ll}
 u_{\theta^1}'' + \sigma\, u_{\theta^1} = \delta_{\theta^1}, & x \in [0,1]  \\
u_{\theta^1}= 0, & x = 0,1,
\end{array}
\right. \quad 
\text{adjoint}: ~~\left\{ 
\begin{array}{ll}
v_{\theta^2}'' + \sigma\, v_{\theta^2} = \delta_{\theta^2}, & x \in [0,1]  \\
v_{\theta^2} = 0, & x = 0,1.
\end{array}
\right.
\]
We use the finite-difference method to solve these equations numerically. 
We represent the function $\sigma(x)$ as a $d$-dimensional vector $\vecsigma$ via a basis function set $\{ \psi_k\}_{k=1}^d$, namely $\sigma(x) \approx \sum_{k=1}^d \sigma_k\, \psi_k(x),~~\text{with}~\vecsigma = (\sigma_1, \dots, \sigma_d)^\top$. In the discrete setting, the above Fr\'echet derivative $\frac{\delta \M}{\delta \vecsigma}$ transfers to a gradient vector with respect to $\vecsigma \in \R^d$, whose $k$-th coordinate is 
\[
\partial_{\sigma_k}\M(\theta;\vecsigma )=\left(\nabla_{\vecsigma}\M(\theta;\vecsigma )\right)_k  = \left\langle\frac{\delta\M(\theta;\sigma )}{\delta \sigma}, \psi_k \right\rangle= \left\langle u_{\theta^1} \cdot v_{\theta^2}, ~\psi_k \right\rangle, \quad \text{for}~ k = 1, \dots, d.
\]
We decompose the domain $[0,1]$ into $d$ equi-spaced cells and set the basis $\{\psi_k\}_{k=1}^d$ to be characteristic functions over these cells. 
The discrete potentials are $\vecsigma_k \equiv \sigma(x_k)$, where $x_k$ is the center of the $k$-th cell for $k = 1, \dots, d$. 
Throughout this section, we set $d=100$.

{\bf Algorithm setup.}
The streamlined gradient flow method, \cref{alternative alg}, alternates between updates for the design particles $\{\theta_i\}_{i=1}^N$ and the reconstruction discrete vector $\vecsigma^* \in \R^{d}$. 
For both design criteria (A- and D-optimal), we choose a ground truth $\vecsigma_\true$, and define the observation data to be $\text{data}(\theta) = \M(\theta; \vecsigma_{\text{true}})+\epsilon$, with the additive noise $\epsilon$ drawn from  $0.1\, \mathcal{N} (0, 1)$. 
We consider two types of choices for the initial design measure $\rho^0$: a uniform distribution on $\Omega = [0,1]^2$ and a warm-start distribution. 
Given $\rho^0$, the initial reconstruction result $\vecsigma^{*,0}$ \eqref{eqn: continuous loss} is obtained by running the inner loop of \cref{brute-force alg} until convergence to a critical point, starting from a random guess drawn from $\vecsigma^0 \sim \vecsigma_{\text{true}} + 2\, U[0,1]^{100}$.
In one simulation, we randomly draw $N$ particles $\{\theta_i \}_{i=1}^N$ to realize the initial design measure $\rho^0$. 
We conduct $10$ independent simulations for each numerical experiment and plot their averaged results in the figures below (with the shading around curve representing the standard deviation of the algorithm runs). 

In the main loop of \cref{alternative alg}, to reduce the computation burden, we utilize an approximate calculation for the velocity term $\dot{\vecsigma}^*$ \eqref{eqn: sigma velocity}.
Recognizing that the distance between $\vecsigma^*$ and the ground-truth $\vecsigma_{\text{true}}$ is small, so that the discrepancy $\M(\theta; \vecsigma^*) - \text{data}(\theta)$ is small,  we drop the second lines of both formulas \eqref{eqn: hess loss} and \eqref{eqn: partial theta}.
By doing so, we avoid the nontrivial calculations of terms $\text{Hess}_{\vecsigma} \M$ and $\nabla_\theta \nabla_{\vecsigma} \M.$


\subsection{A-optimal design} \label{subsec: numerics IA}
First we consider the case of A-optimal design \eqref{eqn: A-optimal_star}. 
We will demonstrate two tests with different initialization choices for the design measure $\rho^0$. 
The ground-truth potential $\sigma_{\text{true}}$ is chosen to be the Gaussian mixture illustrated in the following figure.

\begin{figure}[!htb]
\centering
\includegraphics[scale=0.3]{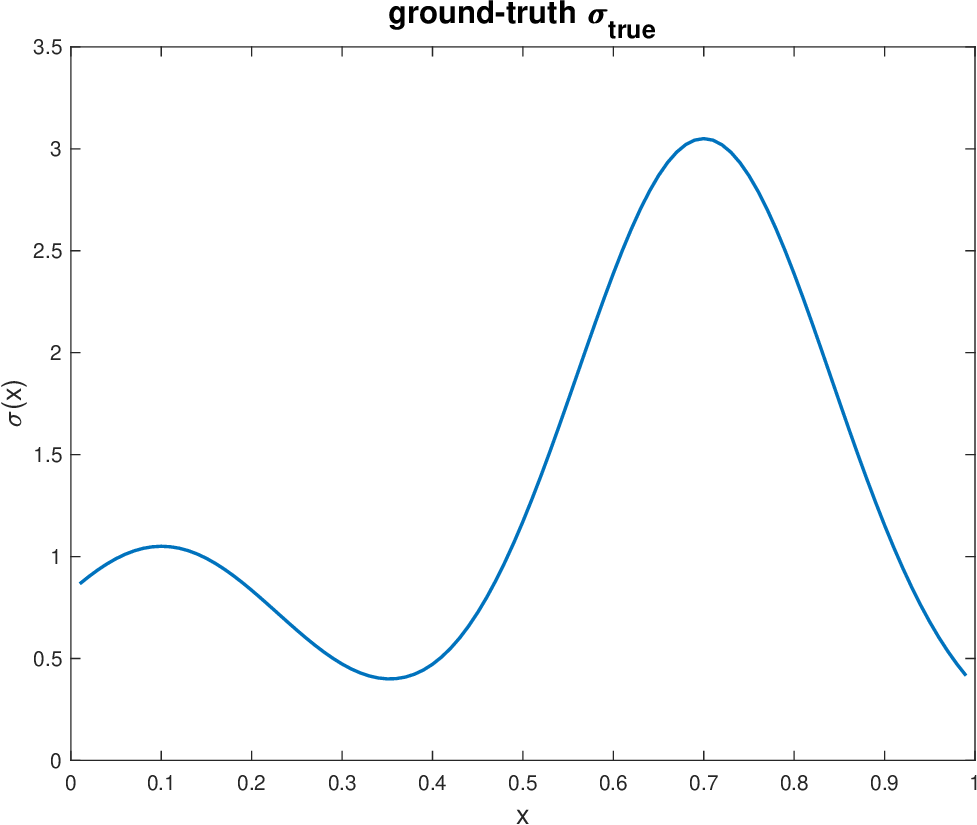}
\caption{Ground-truth $\sigma_\true=\text{exp}\left(-\frac{x-0.1}{0.2} \right)^2+3\,\text{exp}\left(-\frac{x-0.7}{0.2} \right)^2 +0.05$.}
\label{fig: A media}
\end{figure}


In the first test, the initial design measure is the uniform distribution $\rho^0  \sim U ([0,1]^2)$. 
We sample $N = 10,000$ particles uniformly on the design space $\Omega =[0,1]^2$. 
For  the main loop of \cref{alternative alg}, we set $T = 2000$ and  $\Delta t = 10^{-7}$.

\begin{figure}[!htb]
    \centering  \includegraphics[scale=0.4]{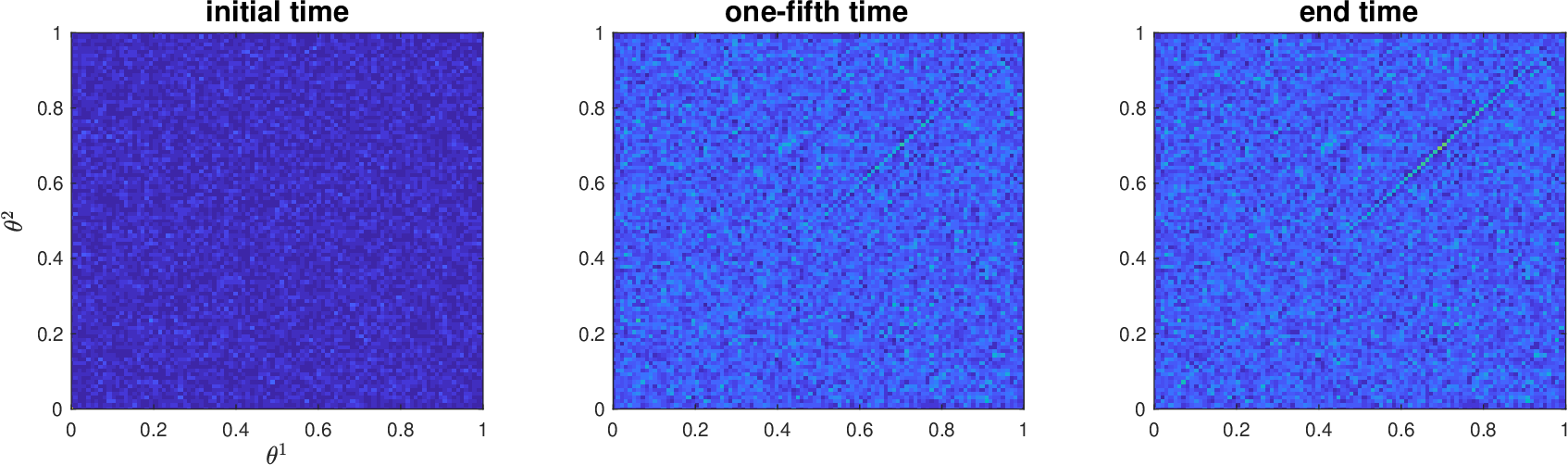}
    \caption{A-optimal design measure at three different points in the algorithm's progress, when  $\rho^0$ is drawn from $U[0,1]^2$.}
    \label{fig: A uniform}
\end{figure}

\begin{figure}[!htb]
\centering   \includegraphics[scale=0.48]{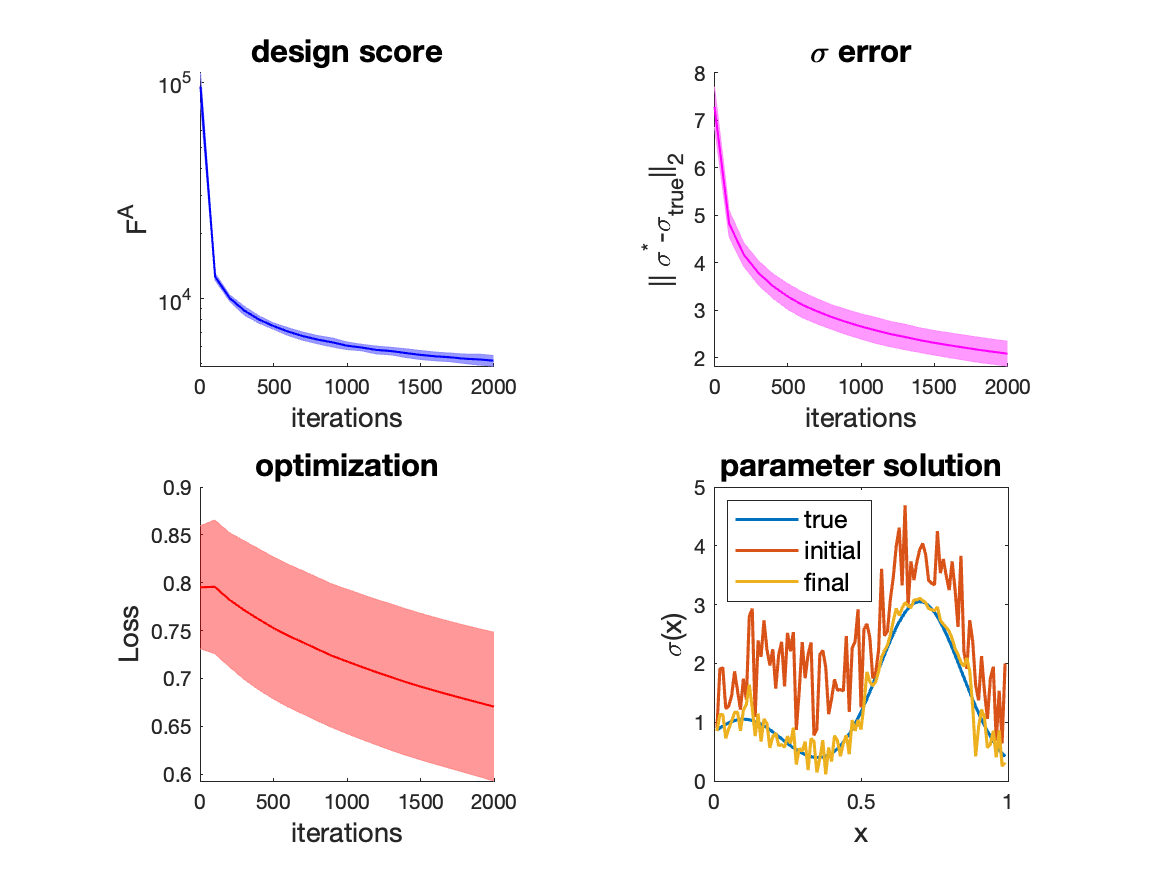}
   \vspace{-1em}
    \caption{Performance evaluation of A-optimal design: uniform distribution initialization.}
    \label{fig: A uniform score}
\end{figure}

We showcase three snapshots of the design measure $\rho^t$ in \cref{fig: A uniform}. 
As the density  $\rho$ evolves, it gradually picks up the diagonal information, indicating that the measurements are the most informative for Schr\"odinger reconstruction when the source $\theta^1$ and detector $\theta^2$ are in close proximity. Note that \cref{fig: A uniform} demonstrate the density result by counting the particles from 10 simulations, where each run returns very similar measure evolution.
In~\cref{fig: A uniform score}, we plot three metrics: A-optimal score $F^A$ \eqref{eqn: A-optimal_star} (in log-scale), $L^2$ reconstruction error $\| \vecsigma^* - \vecsigma_\true \|_2$, and the loss function \eqref{eqn: continuous loss}.
They all decrease with iteration count, suggesting that the adaptive gradient flow \cref{alternative alg} works effectively for A-optimal design as well as the inverse problem. 
In the  bottom right panel of~\cref{fig: A uniform score}, we compare the reconstructed $\vecsigma^*$ at the beginning and end of the algorithm to the ground-truth. 
Note that although the initial reconstruction $\vecsigma^{*,0}$ (red) is jagged and bears little resemblance to $\vecsigma_{\text{true}}$ (blue), the final solution $\vecsigma^{*,T}$ (yellow) closely fits the ground-truth. All the four plots in \cref{fig: A uniform score} also present the averaged metric and parameter values from all simulations.

\vspace{-2em}

In the second test, we use warm-start initialization, setting $\rho^0$ to be close to the optimal design measure obtained from the first test. 
As observed from the output in~\cref{fig: A uniform}, the design measure $\rho$ has pronounced concentration along the diagonal line on the design space $[0,1]^2$, implying that source-detection pairs with $\theta^1 = \theta^2$ contain more valuable information for reconstruction. 
Thus, in the second test, we begin by sampling the bivariate particles $(\theta^1, \theta^2)$ according to $\theta^1 - \theta^2 \sim 0.002\,\mathcal{N}(0,1)$. 
The number of sampled particles is set as $N = 1000$. We pursue the A-optimal design by running \cref{alternative alg} with step size $\Delta t =10^{-6}$ and total number of iterations $T = 1000$. 

\begin{figure}[!htb]
    \centering
 
 \subfloat[\footnotesize{A-optimal design measure at the start and final stages of \cref{alternative alg}, with the initialization $\rho^0$ concentrating around the diagonal stripe.}]{\includegraphics[scale=0.4]{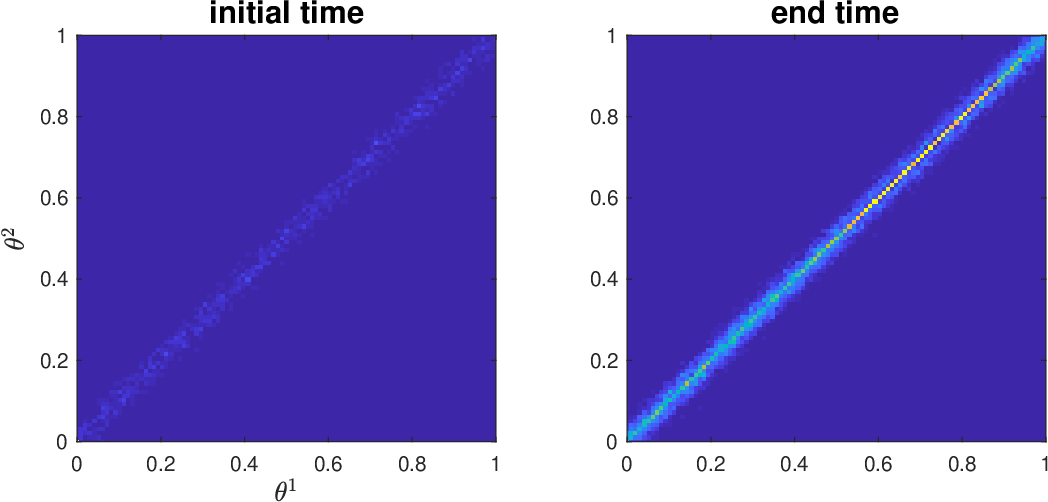}}~~~
 \subfloat[\footnotesize{Particle histogram on the diagonal line at the final stage.}]{\includegraphics[scale=0.25]{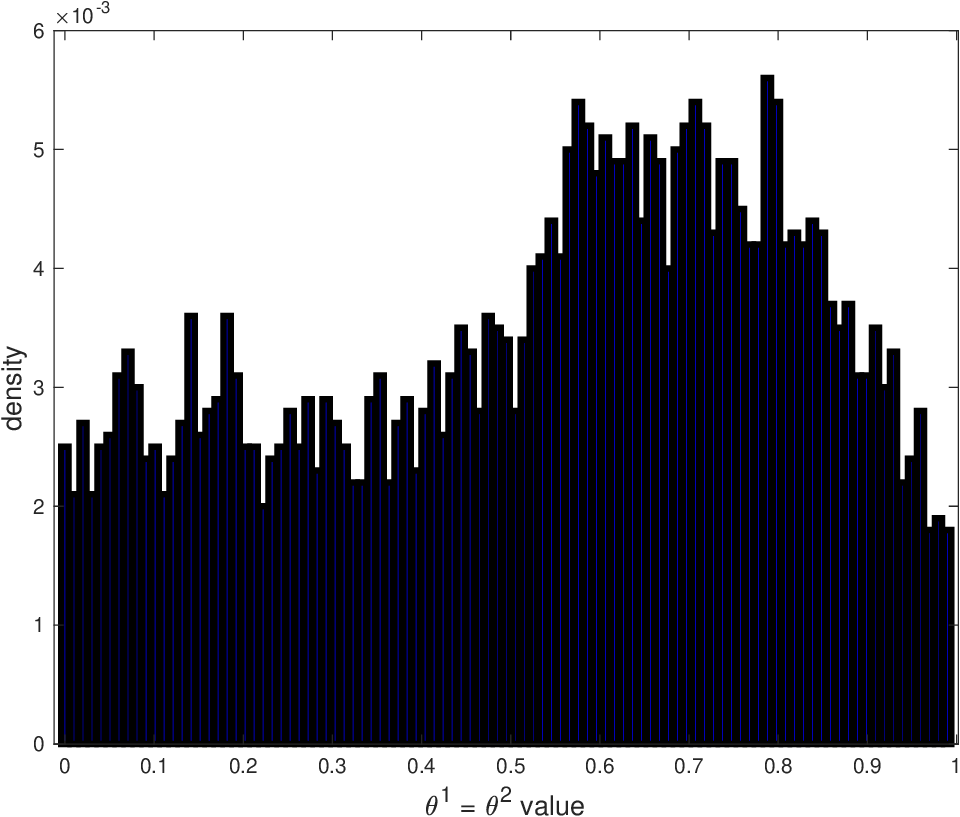}}
 \vspace{-1.5em} 
 \caption{A-optimal warm-start case.}
 \label{fig: A diagonal}
\end{figure}

\begin{figure}[!htb]
\centering
\includegraphics[scale=0.46]{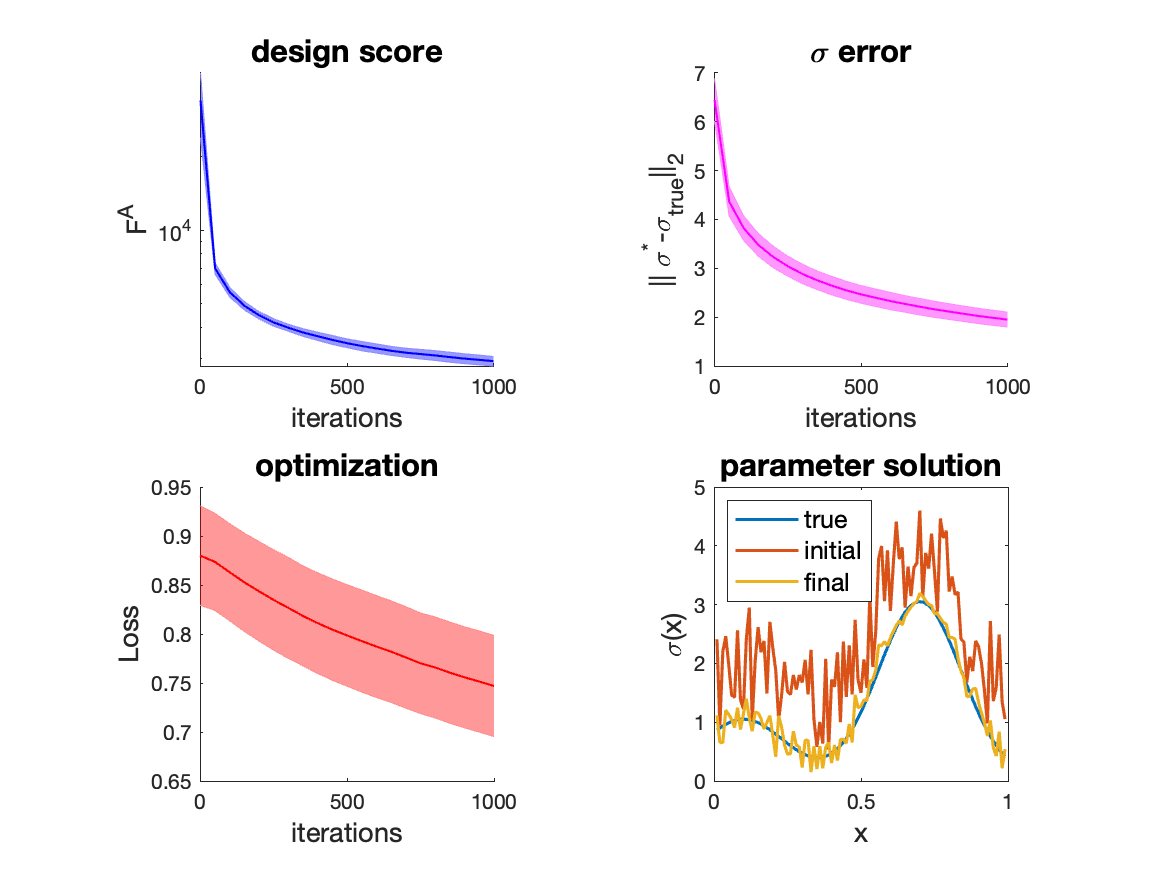}
\vspace{-1em}
    \caption{Performance evaluation of A-optimal design:  warm-start case.} 
    \label{fig: A diagonal score}
\end{figure}

The numerical results are shown in~\cref{fig: A diagonal}. 
In panel (a) the initial measure and final measures from \cref{alternative alg} appear side by side, showing that the algorithm  concentrates the particle even more  along the diagonal. 
Panel (b) displays the histogram of particles constrained on the 1D diagonal line $\{\theta^1 = \theta^2 \}$. 
The constrained distribution roughly follows the ground-truth potential function shape of \cref{fig: A media}. 
It implies that the most informative design is obtained by setting the source and detector at locations where $\vecsigma_\true$ has higher values. 
Performance evaluation of the algorithm is shown in \cref{fig: A diagonal score}.
The A-design score, the reconstruction error of $\vecsigma$, and the loss function all decay, while the bottom-left panel shows the final output density aligning well with  the ground truth $\vecsigma_\true$. 
   \vspace{-1em}
   
\subsection{D-optimal design}
\label{subsec: numerics ID}
\begin{figure}[!htb]
\centering
\includegraphics[scale=0.3]{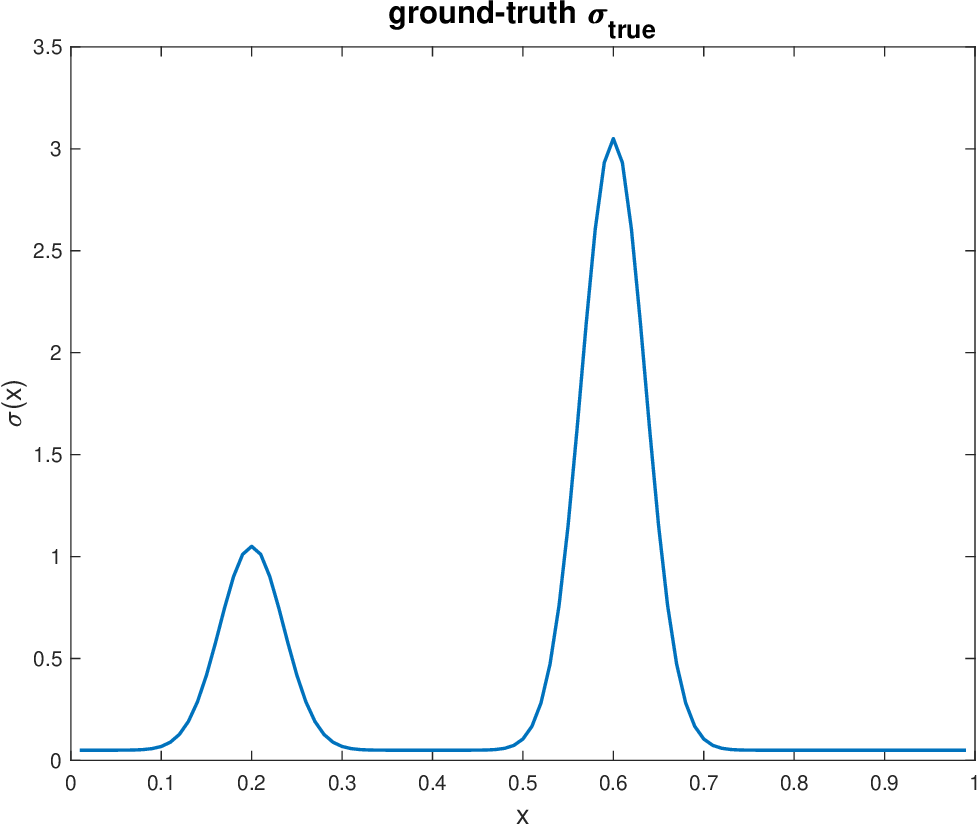}
\vspace{-0.5em}
\caption{Ground-truth $\sigma_\true= \text{exp}\left(-\frac{x-0.2}{0.05} \right)^2+3\,\text{exp}\left(-\frac{x-0.6}{0.05} \right)^2+0.05$.}
\label{fig: D media}
\end{figure}

We now illustrate the numerical result of \cref{alternative alg}  on the D-optimal design criterion \eqref{eqn: D-optimal_star}. 
As for the A-optimal case (\cref{subsec: numerics IA}), we present two testing examples, in which the design measure is initialized with a  uniform distribution and a warm-start, respectively. 
In the D-optimal case, we set the ground-truth potential $\sigma_\true$ to be Gaussian mixture (see \cref{fig: D media}), but with sharper clusters than in the A-optimal setup of \cref{fig: A media}.

\begin{figure}[!htb]
    \centering    \includegraphics[scale=0.35]{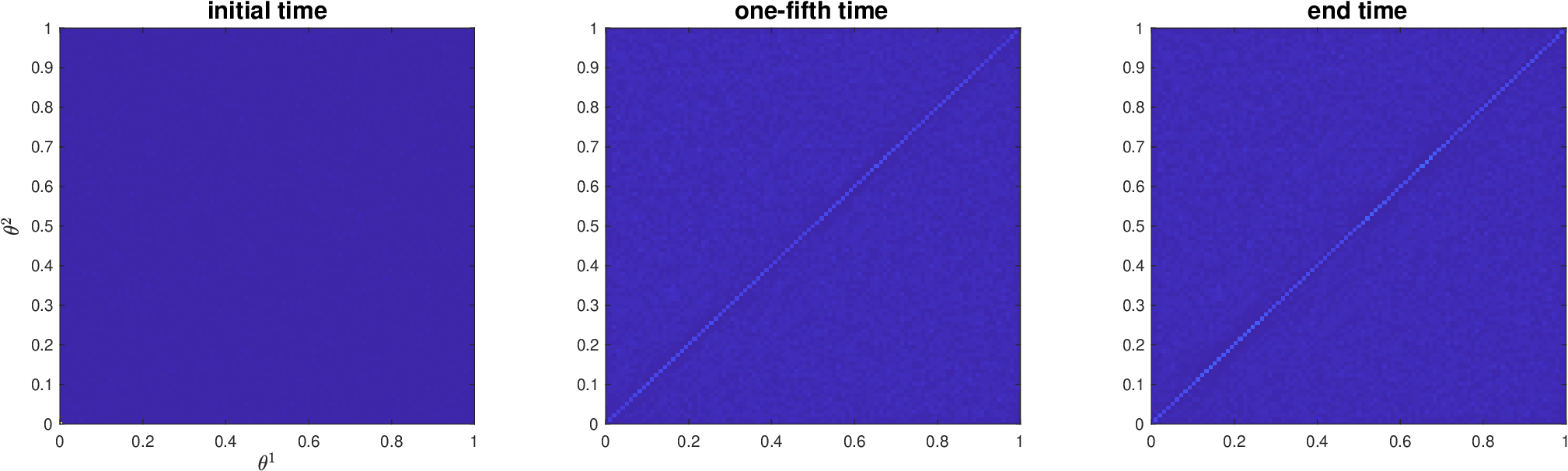}
    \vspace{-1em}
    \caption{D-optimal design measure at three different snapshots when $\rho^0$ is the uniform distribution.}
    \label{fig: D uniform}
\end{figure}

For the uniform distribution initialization $\rho^0  \sim U ([0,1]^2)$,  we draw  $N = 10,000$ samples uniformly from the design space. 
We execute \cref{alternative alg} with  $T = 500$ iterations and step size $\Delta t =10^{-2}$. 
\cref{fig: D uniform} shows the evolution of the design measure $\rho$ as the algorithm proceeds. 
As for the A-optimal design results of \cref{subsec: numerics IA}, the particles gradually gather around the diagonal line, suggesting that coincident source and detection locations provide valuable information under the D-optimal criterion too.

\begin{figure}[H]
\centering \vspace{-1em}
 \includegraphics[scale=0.5]{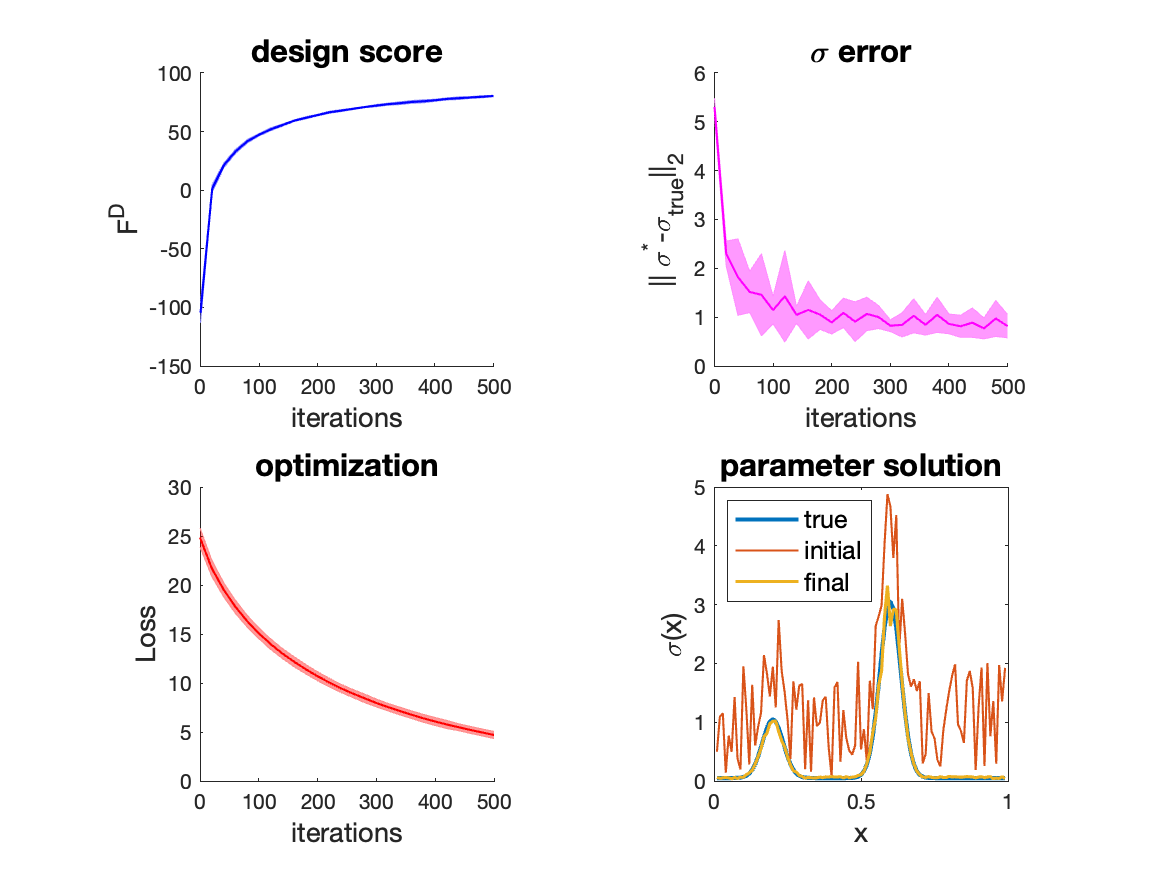}
    \caption{Performance evaluation of D-optimal design, algorithm initialized with uniform distribution.}
    \label{fig: D uniform score}
\end{figure}

Computational performance of \cref{alternative alg} is shown in~\cref{fig: D uniform score}.
The D-optimal design score \eqref{eqn: D-optimal_star} rises to a plateau, while reconstruction errors for $\vecsigma$ and the loss function both decrease. 
It is also noticeable that the $\|\vecsigma-\vecsigma_{\text{true}} \|_2$ in the upper-right panel is more volatile than for the A-optimal design case (\cref{fig: A uniform score}). 
In the bottom-right panel of \cref{fig: D uniform score}, we observe that even though the initial solution $\vecsigma^{*,0}$ (red) is far from the ground-truth $\vecsigma_\true$ (blue), our final output (yellow) of \cref{alternative alg} matches the true potential closely.

Next, we use a warm-start strategy, initializing $\rho^0$ to concentrate along the diagonal line of the design space $\Omega = [0,1]^2$. 
We sample $N = 1000$ particles across the diagonal stripe, exactly as the A-optimal setup described in~\cref{subsec: numerics IA}. 
We proceed to run \cref{alternative alg} for $T=500$ iterations with step size $\Delta t = 10^{-2}$.

\begin{figure}[!htb]
\centering
\subfloat[\footnotesize{D-optimal design measure at the initial and final stages of \cref{alternative alg}, when the initial $\rho^0$ is set around the diagonal stripe.}]
{ \includegraphics[scale=0.37]{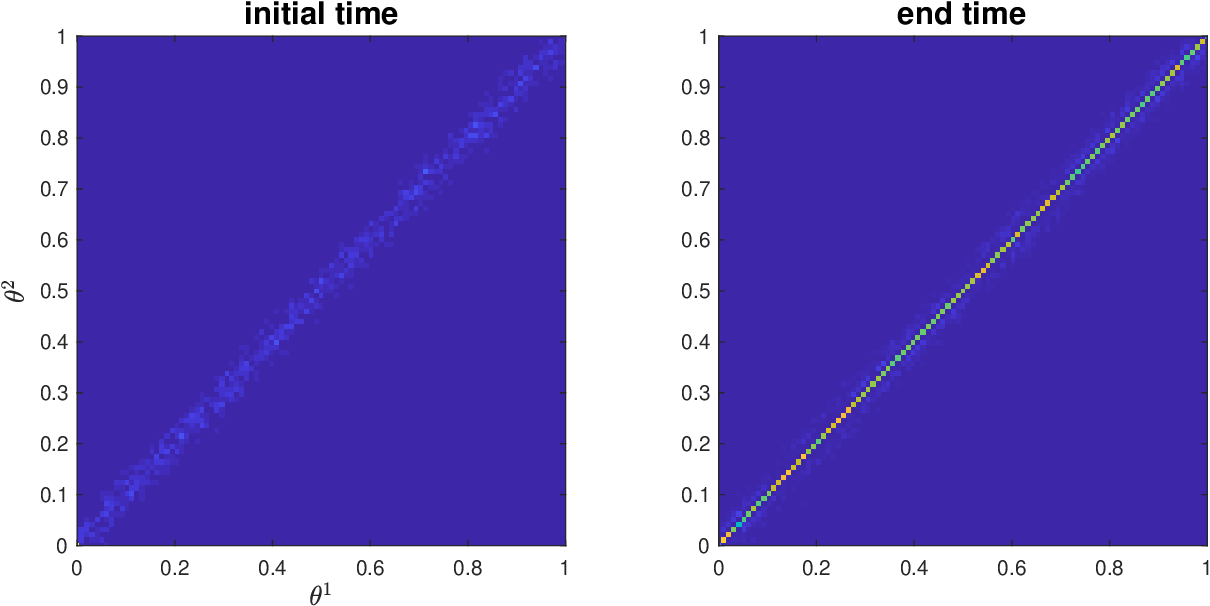}}~~~~
\subfloat[\footnotesize{Particle histogram on the diagonal line at the final stage.}]{   \includegraphics[scale=0.28]{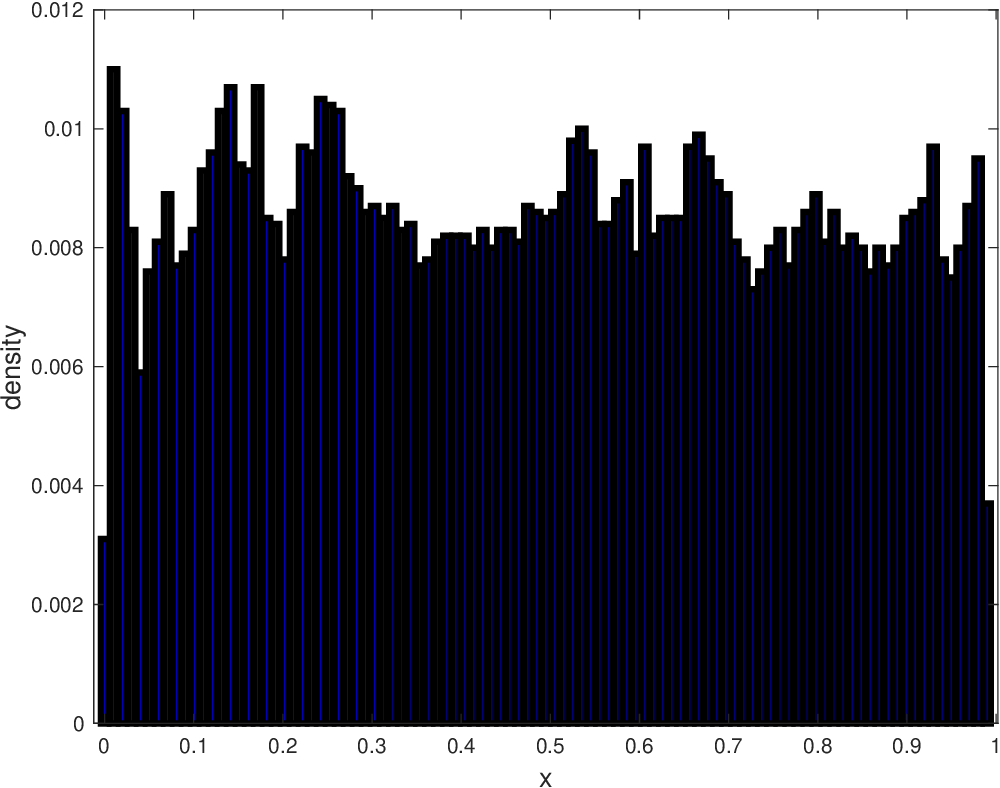}}
    \caption{D-optimal warm-start case.}
    \label{fig: D diag}
\end{figure}

\begin{figure}[!htb]
\centering
\includegraphics[scale=0.5]{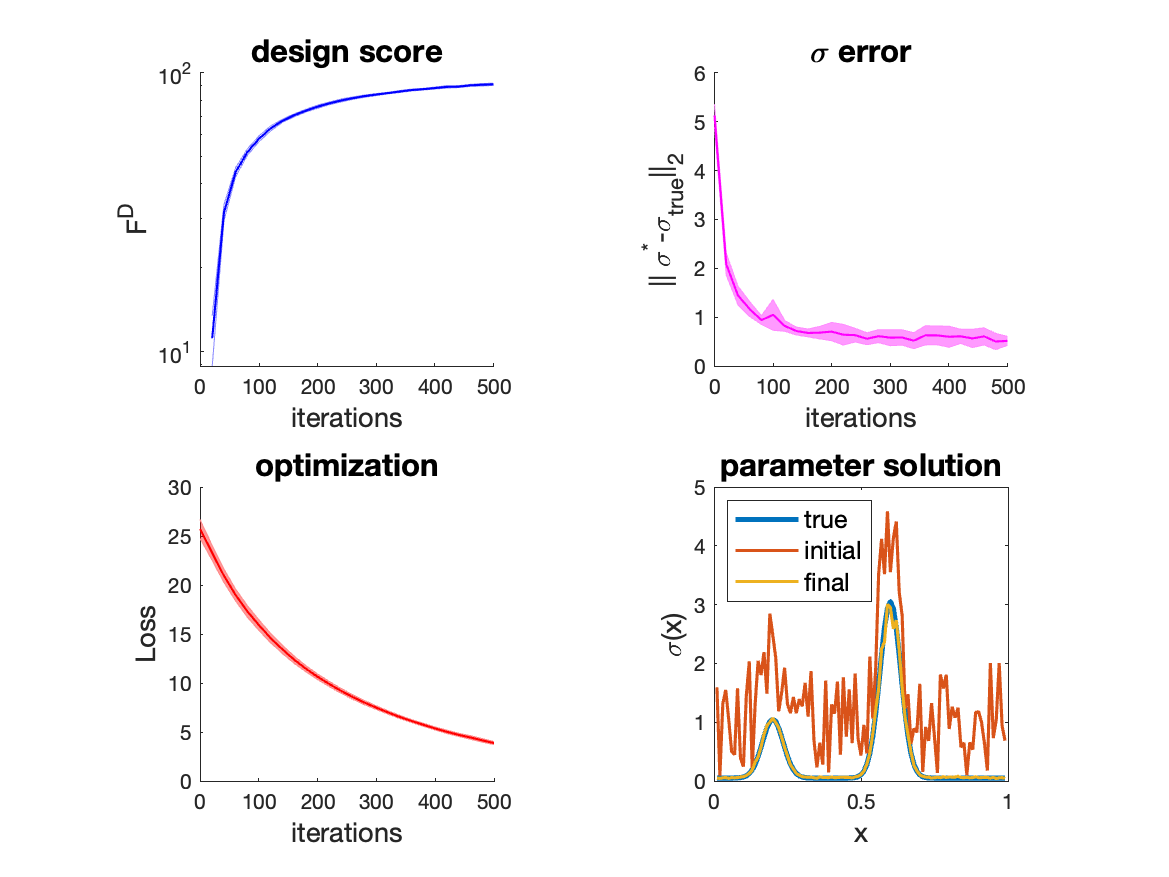}
    \caption{Performance evaluation of D-optimal design: warm-start case.}
    \label{fig: D diag score}
\end{figure}

Part (a) of~\cref{fig: D diag} shows initial and final design measure $\rho$, while part (b) shows the histogram of particles constrained on the 1D diagonal line $\{\theta^1 = \theta^2 \}$. 
Again, we see the concentration along the diagonal, but unlike for the A-optimal design criterion, the D-optimal design shows no particularly favored locations along the diagonal, even though the true potential $\vecsigma_{\text{true}}$ (\cref{fig: D media}) has widely varying  values across this interval. This suggests the D-optimal criterion is not sensitive to the influence of potential value distribution.   
In~\cref{fig: D diag score}, we  show computational performance of the algorithm.
The results are roughly similar to the random initialization case (\cref{fig: D uniform score}), although the top-right panel shows that the design score starts from a value much close to optimality, as expected.

\clearpage
\section{Conclusions and future directions}

Efficient data acquisition is crucial for learning physical systems. 
By computing an optimal probability distribution that quantifies the importance of data, we can weight the data accordingly to best extract reconstruction of unknown parameters.
Novel gradient flow and particle simulation methods can calculate experimental designs over probability measure space $\pr_2$. 
Since most physical models are highly nonlinear, the design measure $\rho$ and the vector $\vecsigma^*$ of  parameters to be recovered are co-dependent. 
We proposed a bilevel optimization scheme that adaptively solves for $\rho$ and $\vecsigma^*.$ 
The difference between our two algorithmic approaches --- brute-force and streamlined --- lies in the update of unknowns $\vecsigma^*$ for a given measure $\rho$. 
The streamlined formulation employs the one-step update for $\vecsigma^*$, reducing computational cost. 
We applied these methods to the Lorenz dynamical system and a Schr\"{o}dinger model. 
Our algorithms returned optimal design measures and revealed interesting physics phenomena.

Despite promising empirical performance, the current work still has several limitations described below, which will be the subjects of future studies.
\begin{enumerate}
\item {\bf Convergence theory.} The convergence properties of the adaptive gradient flow method (\cref{alternative alg}) need to be understood better. 
Open problems include the dependence of these properties on the sample size $N$ and number of iterations $T$. 
The dependence between the design distribution $\rho$ and the reconstruction $\vecsigma^*[\rho]$ prompts analytical challenges. 
Is it possible to examine the convergence of \cref{alternative alg} via studying the coupled updating system of $\rho$ and $\vecsigma^*$? 
Does evolution of the design measure $\rho$ lead to the reconstruction $\vecsigma^*[\rho]$ converging to the ground-truth $\vecsigma_\true?$

\item {\bf Efficiency improvement.} 
In our tests, hyperparameters such as $T$ and $\Delta t$ are tuned manually, but these should be chosen automatically and adaptively in practice. 
Moreover, when the sample particle $\theta$ gets updated, the algorithm has no further use for the measurement $\text{data}(\theta)$.
Since these observations are costly, can we develop a scheme to make better use of the measurements and recycle them in certain way?
\end{enumerate}
They are all very interesting questions that need to be answered to fully understand the computational and sample complexity of the algorithm and its effectiveness. We leave these questions for future research.

\section*{Acknowledgements}
R.J. and Q.L. acknowledge support from NSF-DMS-2308440. R.J. is supported by NSF-DMS-2023239.  
S.W. acknowledges support from NSF DMS-2023239 and CCF-2224213.
\appendix 
\section{Computation details for the Lorenz model}
\label{sec: appendix}


We  provide some  computational details concerning the Lorenz 63 model. 

\begin{figure}[!htb]
\includegraphics[scale=0.3]{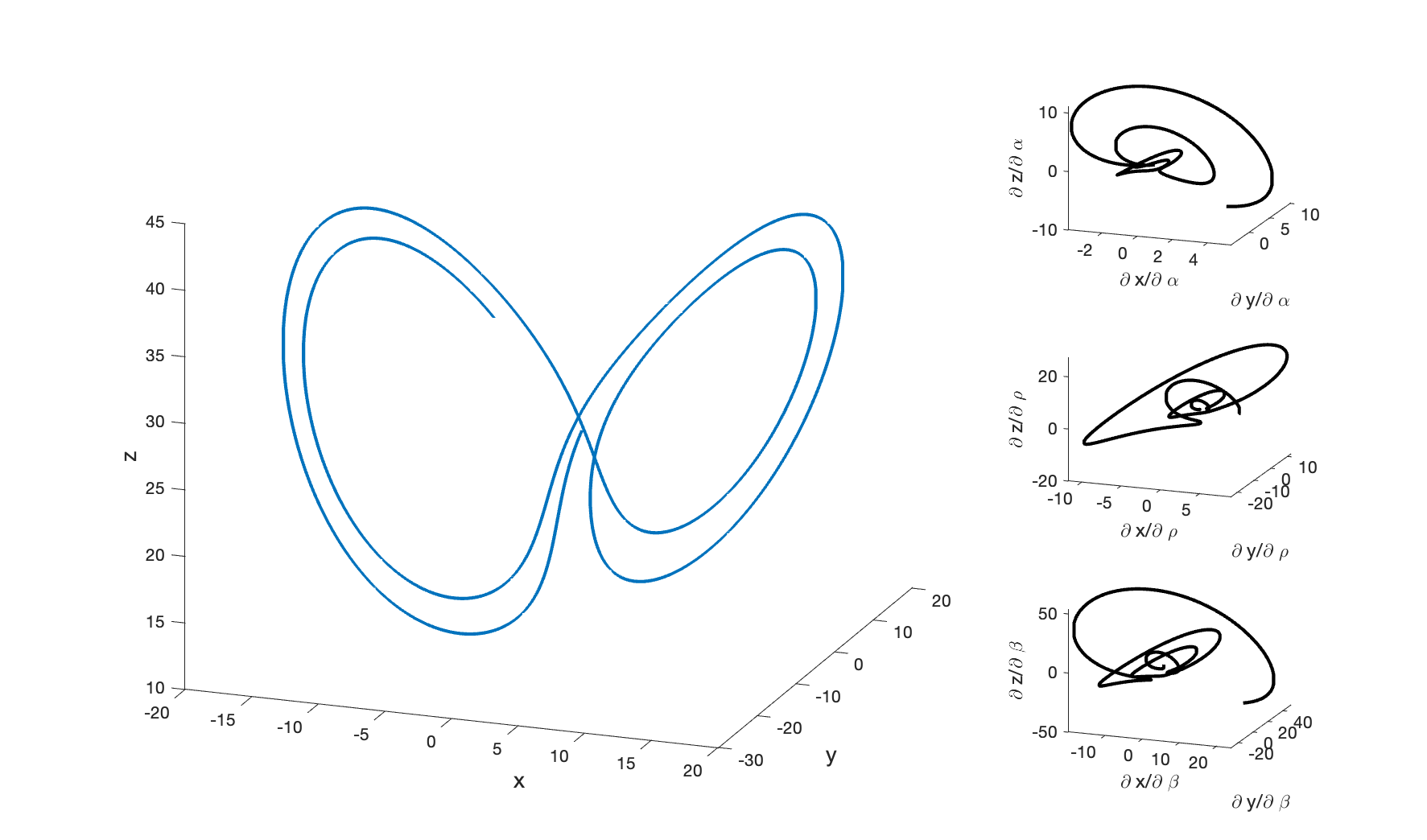}
    \caption{Lorenz 63 model \eqref{eqn: L63}  (left) and its gradient terms w.r.t. $\vecsigma$ \eqref{eqn: L63 derivatives} (right).}
    \label{fig: L63 model}
\end{figure}

The experimental design objectives \eqref{eqn: A,D-optimal_star} rely on the gradient term $\nabla_{\vecsigma} \M \in \R^3$ with respect to the model parameters $\vecsigma = (\alpha, \gamma, \beta)$.  
This gradient can be computed explicitly as one of the rows of the following $3 \times 3$ matrix depending on the observable $x, y, z:$
\begin{equation}
\label{eqn: L63 Mgrad}
J := 
\left[\begin{array}{ccc}
\vspace{1mm}
\frac{\partial x}{\partial \alpha}& \frac{\partial x}{\partial \gamma} & \frac{\partial x}{\partial \beta}\\
\vspace{1mm}
\frac{\partial y}{\partial \alpha} &\frac{\partial y}{\partial \gamma} &\frac{\partial y}{\partial \beta} \\
\frac{\partial z}{\partial \alpha} & \frac{\partial z}{\partial \gamma} &\frac{\partial z}{\partial \beta}
\end{array}
\right]\,.
\end{equation}
Differentiating the Lorenz model~\eqref{eqn: L63} on both sides with respect to $\vecsigma = (\alpha, \gamma, \beta)$, we obtain
\begin{equation}
\label{eqn: L63 derivatives}
\frac{\rd}{\rd \tau}J
= 
\left[\begin{array}{ccc}
  -\alpha & \alpha & 0   \\
  \gamma-z  & -1&-x\\
  y&x&-\beta
\end{array}
\right]\cdot J
 +
\left[
\begin{array}{ccc}
y-x & 0 & 0\\
0 &x & 0\\
0 & 0 & -z\\
\end{array}
\right]\,.
\end{equation}
Since the initialization for the state variables are fixed, the initial value for these derivatives are all 0, meaning that $J(\tau=0) = {\bf 0}$.
Solving this dynamical system gives us the evaluation of $J$ for all times $t$. 
In \cref{fig: L63 model} we plot the trajectory of Lorenz 63 model with three state variables and their corresponding gradients w.r.t. $\vecsigma$, evaluated at the ground-truth $\vecsigma_\true$. 
A similar strategy applies to computing the second-order term $\text{Hess}_{\vecsigma} \M$ in \eqref{eqn: hess loss}. 

In the implementation of our particle gradient flow algorithms, we also need the gradients w.r.t. the design variable $\theta$, including the terms $\nabla_\theta \M$ and $\nabla_\theta \nabla_{\vecsigma} \M$. Once the observable $c\in \{x, y, z\}$ is determined, the differentiation of $\theta = (c, \tau)$ is simplified as the time derivative $\frac{\dd}{\dd \tau} (\cdot)$. Hence $\nabla_\theta \M$ and $\nabla_\theta \nabla_{\vecsigma} \M$ can be extracted respectively from the left-hand sides of the dynamical systems \eqref{eqn: L63} and \eqref{eqn: L63 derivatives}.

\bibliographystyle{plain}
\bibliography{references}

\end{document}